\documentclass{article}
\usepackage{amssymb}
\usepackage{amsmath}
\usepackage{amsfonts}
\usepackage{amscd}
\usepackage[dvips]{graphicx}
\usepackage{amssymb,amsfonts,amsxtra, mathrsfs,placeins,graphicx,verbatim}
\usepackage[all]{xy}

\xyoption{line}

\newtheorem{lemm}{Lemme}

\newtheorem{corollary}{Corollary}
\newtheorem{theorem}{Theorem}
\newtheorem{definition}{Definition}
\newtheorem{proposition}{Proposition}

\begin{document}

\title{$n$-ary  associative algebras, cohomology, free algebras and coalgebras}
\author{Nicolas GOZE, Elisabeth REMM \\
 nicolas.goze@uha.fr, elisabeth.remm@uha.fr\\
adresse : LMIA.\ UHA\\
4 rue des Fr\`{e}res Lumi\`{e}re\\
F. 68093 Mulhouse Cedex}
\maketitle

\begin{abstract}
When $n$ is odd, a cohomology of type Hochschild for $n$-ary partially
associative algebras has been defined in Gnedbaye's thesis. Unfortunately,
the cohomology definition is not valid when $n$ is even. This fact is found
again in the computations of the $n$-ary partially associative free algebra.

In this work, we define in a first time two approachs of an Hochschild
cohomology for $n$-ary partially associative algebras. First by reducing the
space of cochains, secondly by using a graded version. Next we compute the
free $n$-ary algebra, giving a basis of this algebra. At last we extend the notion of coalgebras to $n$-ary algebras.
\end{abstract}

All algebraic objects will be considered over a commutative field $\mathbb{K}
$ of characteristic zero.

\section{Relations between $n$-ary partially associative algebras and
Gerstenhaber products}

\subsection{\protect\bigskip Definition}

Let $V$ be a $\mathbb{K}$-vector space and consider 
\begin{equation*}
C^{k}(V)=Hom_{\mathbb{K}}(V^{\otimes k},V),
\end{equation*}%
for any natural number $k$. By defintion a $n$-ary partially associative
algebra is a pair $(V,\mu )$ where $\ V$ is a $\mathbb{K}$-vector space and $%
\mu $ a linear map $\mu :V^{\otimes n}\rightarrow V$ satisfying%
\begin{equation} 
\label{equa}
\underset{i=1}{\overset{n}{\sum }}(-1)^{(i-1)(n-1)}\mu (X_{1},\cdots ,\mu
(X_{i},\cdots ,X_{i+n-1}),X_{i+n-1}),\cdots ,X_{2n-1})=0.
\end{equation}%
When 
\begin{equation*}
\mu (X_{1},\cdots ,\mu (X_{i},\cdots ,X_{i+n-1}),\cdots ,X_{2n-1})=\mu
(X_{1},\cdots ,\mu (X_{j},\cdots ,X_{j+n-1}),\cdots ,X_{2n-1})
\end{equation*}%
for any $i,j\in \{ 1,\cdots ,p\},$ the algebra $(V,\mu )$ is
totally associative.

\noindent 

\subsection{Gerstenhaber products $\bullet _{n,n}$}

These products have been proposed by Gerstenhaber  in the study of  spaces
 of Hochschild cohomology of an associative algebra. We recall this quickly
in order to use the practical notations which appear in the work of
Gerstenhaber.

The Gerstnhaber product of $f\in C^{n}(V)$ and $g\in C^{m}(V)$ is the
element $f\bullet _{n,m}g\in C^{n+m}(V)$ defined by

\begin{equation*}
f\bullet _{n,m}g(X_{1}\otimes \cdots \otimes
X_{n+m-1})=\sum_{i=1}^{n}(-1)^{(i-1)(m-1)}f(X_{1}\otimes \cdots \otimes
g(X_{i}\otimes \cdots \otimes X_{i+m-1})\otimes \cdots \otimes X_{n+m-1}).
\end{equation*}%
These Gerstenhaber products satisfy pre-Lie Identity [see Ge] that is:

\begin{equation*}
(f\bullet_{n,m}g)\bullet_{n+m-1,p}h -f\bullet_{n,m+p-1}(g \bullet_{m,p}h) =
(-1)^{(m-1)(p-1)}\left((f\bullet_{n,p}h)\bullet_{n+p-1,m}g -
f\bullet_{n,n+p-1}(h\bullet_{p,m}g\right),
\end{equation*}
for any $f\in C^{n}(V),\quad g\in C^{m}(V)$ and $h\in C^{p}(V).$

\medskip

\noindent \textbf{Notations.} We denote the products of Gerstenhaber by $%
\bullet _{n,k}$. When there is no confusion, we denote these products simply
by $\bullet $. Moreover, the symbol $\circ $ refers to the ordinary
composition of applications.

\begin{definition}
We call $n$-ary algebra associated to $\bullet _{n,n}$ any $\mathbb{K}$%
-vector space $V$ with an application $\mu \in C^{n}(V)$ satisfying: 
\begin{equation*}
\mu \bullet _{n,n}\mu =0.
\end{equation*}
We denote it by $(V,\bullet _{n,n})$ or $(V,\bullet _{n,n},\mu)$ if we need
to specify the multiplication $\mu $.
\end{definition}

\noindent Then we have:

\begin{eqnarray} \label{nary}
\mu \bullet _{n,n}\mu (X_{1},\cdots ,X_{2n-1})=\underset{i=1}{\overset{n}{%
\sum }}(-1)^{(i-1)(n-1)}\mu (X_{1},\cdots ,\mu (X_{i},\cdots
,X_{i+n-1}),\cdots ,X_{2n-1})=0.
\end{eqnarray}%
These algebras correspond to partially associative algebras.

\bigskip 

\textbf{Remarks.} We will study identities which are deduced from the
definition of the product $\bullet _{n,n}.$

1. For $n=1$, Identity (\ref{nary}) reduces to: 
\begin{equation*}
\mu \bullet _{1,1}\mu (X_{1})=\mu (\mu (X_{1}))=0
\end{equation*}%
so $\mu \circ \mu =0.$

\bigskip

2. \noindent For $n=2$ we get:

\begin{equation*}
\mu \bullet _{2,2}\mu (X_{1},X_{2},X_{3})=\mu (\mu (X_{1},X_{2}),X_{3})-\mu
(X_{1},\mu (X_{2},X_{3}))=0
\end{equation*}%
and $(V,\bullet _{2,2},\mu )$ is an associative algebra.

3. For $n>2$, the algebra $(V,\bullet _{n,n},\mu )$ corresponds to a
partially associative algebra (with operation in degree $0$) studied by
Gnedbaye (see [Gn]). 

4. A $n$-ary algebra $(V,\bullet _{n,n},\lambda )$ can not be deduced from
an associative algebra $(V,\bullet _{2,2},\mu )$ by composition that is $%
\lambda $ can not be equal to 
\begin{equation*}
\sum_{k=1}^{l}a_{k}\underset{}{(Id_{p_{n-1}}\otimes \mu \otimes
Id_{n-2-p_{n-1}})\circ \cdots \circ (Id_{p_{2}}\otimes \mu \otimes
Id_{n-2-p_{2}})\circ }(Id_{p_{1}}\otimes \mu \otimes Id_{n-p_{1}-2}).
\end{equation*}

5. A $n$-ary Lie algebra [Fi] is defined by a skew-symmetric product $\mu $
satisfying the generalized Jacobi Identity 
\begin{equation*}
\sum_{\sigma \in Sh_{n,n-1}}(-1)^{\varepsilon (\sigma )}\mu (\mu (X_{\sigma
(1)},\cdots ,X_{\sigma (n)}),X_{\sigma (n+1)},\cdots ,X_{\sigma (2n-1)})=0
\end{equation*}%
where $Sh_{n,n-1}$ is the set of $(n,n-1)$-shuffles. If $\lambda $ is a $n$%
-ary multiplication satisfying $\lambda \bullet _{n,n}\lambda =0$, the
product $\mu $ defined by

\begin{equation*}
\mu (X_{1},\cdots ,X_{n})=\sum_{\sigma \in S_n} (-1)^{\varepsilon (\sigma)}
\lambda (X_{\sigma(1)},\cdots, X_{\sigma(n)})
\end{equation*}%
is a $n$-ary Lie algebra product.

\medskip

\begin{lemm}
Let $(V,\bullet _{n,n},\mu )$ be a $n$-ary algebra with product $\bullet
_{n,n}$. If $n$ is even

\begin{equation*}
(\varphi \bullet \mu )\bullet \mu =0,
\end{equation*}%
for any $\varphi \in C^{k}(V).$
\end{lemm}

\noindent \textit{Proof.} It follows from pre-Lie identity that

\begin{equation*}
(\varphi \bullet \mu )\bullet \mu -\varphi \bullet (\mu \bullet \mu
)=(-1)^{(n-1)(n-1)}[(\varphi \bullet \mu )\bullet \mu -\varphi \bullet (\mu
\bullet \mu )].
\end{equation*}%
But $\mu \bullet \mu =0$. Then, as $n$ is even, we obtain

\begin{equation*}
(\varphi \bullet \mu )\bullet \mu =-(\varphi \bullet \mu )\bullet \mu
\end{equation*}%
and, using the fact that $char(\mathbb{K) }=0$, this equation reduces to $%
(\varphi \bullet \mu )\bullet \mu =0.$

\medskip

\noindent Observe that for odd $n$, pre-Lie identity is trivial. In this
case we compute $(\varphi \bullet \mu )\bullet \mu .$
Let $\theta _{k}(\mu )$ be the map $V^{\otimes
(2n+k-2)}\longrightarrow V^{\otimes k}$ defined by

\begin{equation*}
\theta _{k}(\mu )=\underset{%
\begin{array}{c}
{\small 0\leq p\leq k-2} \\ 
{\small 0\leq q\leq k-2-p}%
\end{array}%
}{\sum }Id_{p}\otimes \mu \otimes Id_{q}\otimes \mu \otimes Id_{k-p-q-2}
\end{equation*}
where $Id_{0}$ means no operation, for example $Id_{0}\otimes \mu \otimes
Id_{k-2}\otimes \mu \otimes Id_{0}$ is just $\mu \otimes Id_{k-2}\otimes \mu
.$

\begin{lemm}
If $n$ is odd, then for any $\varphi \in C^k(V)$

\begin{equation*}
(\varphi \bullet \mu )\bullet \mu =2\varphi \circ \theta _k(\mu)
\end{equation*}
where $\circ $ is the ordinary composition. In particular $(\varphi \bullet
\mu )\bullet \mu =0$ if and only if $Im\theta _k(\mu )\in Ker\varphi.$
\end{lemm}

\noindent \textit{Proof. } The product $\bullet $ satisfies 
\begin{equation*}
(\varphi \bullet \mu )\bullet \mu (X_{1},\cdots
,X_{k+2n-2})=\sum_{i=1}^{k+n-1}(\varphi \bullet \mu )(X_{1},\cdots ,\mu
(X_{i},\cdots ,X_{i+n-1}),\cdots ,X_{k+2n-2}).
\end{equation*}

\noindent Terms of the right-hand side are of two kinds. The first
corresponds of elements which can be rewritten as 
\begin{equation*}
\varphi (X_{1},\cdots ,X_{p-1},A,X_{2n+p-1},\cdots ,X_{2n+k-2})
\end{equation*}%
with $1\leq p\leq k$ and 
\begin{equation*}
\begin{array}{lll}
\medskip A & = & \mu (\mu (X_{p},\cdots ,X_{n+p-1}),X_{n+p},\cdots
,X_{2n+p-2})+\mu (X_{p},\mu (X_{p+1},\cdots ,X_{n+p}),\cdots ,X_{2n+p-2}) \\ 
\medskip &  & +\cdots +\mu (X_{p},\cdots ,X_{n+p-2},\mu (X_{n+p-1},\cdots
,X_{2n+p-2})) \\ 
\medskip & = & \mu \bullet \mu (X_{p},\cdots ,X_{2n+p-2}) \\ 
\medskip & = & 0.%
\end{array}%
\end{equation*}%
The second type corresponds to elements 
\begin{equation*}
\varphi (X_{1},\cdots ,X_{p-1},\mu (X_{p},\cdots X_{n+p-1}),X_{n+p},\cdots
,\mu (X_{q},\cdots ,X_{n+q-1}),\cdots ,X_{k+2n-2})
\end{equation*}
with $1\leq p\leq q-n\leq k-1.$ Then

\begin{equation*}
\begin{array}{l}
\medskip (\varphi \bullet \mu )\bullet \mu (X_{1},\cdots ,X_{k+2n-2})= \\ 
\medskip 2\underset{1\leq p\leq q-n\leq k-1}{\sum }\varphi (X_{1},\cdots
,X_{p-1},\mu (X_{p},\cdots X_{n+p-1}),X_{n+p},\cdots ,\mu (X_{q},\cdots
,X_{n+q-1}),\cdots ,X_{k+2n-2}) \\ 
\medskip =2\varphi (\theta _{k}(\mu ))(X_{1},\cdots ,X_{k+2n-2}).%
\end{array}%
\end{equation*}

\section{Cohomology of partially associative algebras $(V, \bullet_{n,n})$}

Recall that if $n=2$, the Hochschild cohomology of an associative algebra
with multiplication $\mu$ is defined from the coboundary operator 
\begin{equation*}
\begin{array}{l}
\medskip \delta^k:C^k(V)\longrightarrow C^{k+1}(V) \\ 
\delta^k\varphi=(-1)^{k-1}\mu\bullet_{2,k}\varphi-\varphi\bullet_{k,2}\mu.%
\end{array}%
\end{equation*}
Consider a $n$-ary algebra with a multiplication $\mu$ of type $%
\bullet_{n,n} $.

\subsection{First case: $n$ is even}

Let $\varphi $ be in $C^{k}(V).$ The applications $\mu \bullet \varphi $ and 
$\varphi \bullet \mu $ are in $C^{k+n-1}(V)$. We define, for any $i\in
\left\{ 0,\cdots ,n-1\right\} $, the linear map: 
\begin{equation*}
\delta _{i}^{k}:C^{i+k(n-1)}(V)\longrightarrow C^{i+(k+1)(n-1)}(V)
\end{equation*}%
by 
\begin{equation*}
\delta _{i}^{k}(\varphi )=(-1)^{k-1}\mu \bullet _{n,n+k-1}\varphi -\varphi
\bullet _{n+k-1,n}\mu .
\end{equation*}

\begin{theorem}
The maps $\delta _{i}^{k}$ satisfy 
\begin{equation*}
\delta _{i}^{k+1}\circ \delta _{i}^{k}=0,
\end{equation*}%
for any $i=0,1,\cdots ,n-2.$
\end{theorem}

\noindent \textit{Proof.} Consider 
\begin{equation*}
\begin{array}{lll}
(\delta _{i}^{k+1}\circ \delta _{i}^{k})(\varphi ) & = & (-1)^{k}\mu \bullet
((-1)^{k-1}\mu \bullet \varphi -\varphi \bullet \mu )-((-1)^{k-1}\mu \bullet
\varphi -\varphi \bullet \mu )\bullet \mu \\ 
& = & -\mu \bullet (\mu \bullet \varphi )+(-1)^{k+1}\mu \bullet (\varphi
\bullet \mu )+(-1)^{k}(\mu \bullet \varphi )\bullet \mu +(\varphi \bullet
\mu )\bullet \mu.%
\end{array}%
\end{equation*}%
The pre-Lie identity implies that 
\begin{equation*}
(\mu \bullet \mu )\bullet \varphi -\mu \bullet (\mu \bullet \varphi
)=(-1)^{(n-1)(k-1)}((\mu \bullet \varphi )\bullet \mu -\mu \bullet (\varphi
\bullet \mu )).
\end{equation*}%
Since $\mu \bullet \mu =0$ and $n$ is even, we obtain 
\begin{equation*}
-\mu \bullet (\mu \bullet \varphi )=(-1)^{k-1}(\mu \bullet \varphi )\bullet
\mu +(-1)^{k}\mu \bullet (\varphi \bullet \mu )
\end{equation*}%
so 
\begin{equation*}
\begin{array}{lll}
\medskip (\delta _{i}^{k+1}\circ \delta _{i}^{k})(\varphi ) & = & 
((-1)^{k-1}+(-1)^{k})(\mu \bullet \varphi )\bullet \mu \\ 
\medskip &  & +((-1)^{k}+(-1)^{k+1})\mu \bullet (\varphi \bullet \mu
)+(\varphi \bullet \mu )\bullet \mu \\ 
\medskip & = & (\varphi \bullet \mu )\bullet \mu.%
\end{array}%
\end{equation*}%
But, as $n$ is even, Lemma 1 implies that $(\varphi \bullet \mu )\bullet \mu
=0$. So $\delta _{i}^{k+1}\circ \delta _{i}^{k}=0$ and we have the following
families of complexes 
\begin{equation*}
\begin{array}{l}
\medskip \mathcal{C}^{0}(V)\overset{\delta _{0}^{0}}{\longrightarrow }%
\mathcal{C}^{n-1}(V)\overset{\delta _{0}^{1}}{\longrightarrow }\mathcal{C}%
^{2(n-1)}(V)\longrightarrow \cdots \longrightarrow \mathcal{C}^{k(n-1)}(V)%
\overset{\delta _{0}^{k}}{\longrightarrow }\mathcal{C}^{(k+1)(n-1)}(V)%
\longrightarrow \cdots \\ 
\medskip \mathcal{C}^{1}(V)\overset{\delta _{1}^{0}}{\longrightarrow }%
\mathcal{C}^{n-1+1}(V)\overset{\delta _{1}^{1}}{\longrightarrow }\mathcal{C}%
^{2(n-1)+1}(V)\longrightarrow \cdots \longrightarrow \mathcal{C}%
^{k(n-1)+1}(V)\overset{\delta _{1}^{k}}{\longrightarrow }\mathcal{C}%
^{(k+1)(n-1)+1}(V)\longrightarrow \cdots \\ 
\vdots \\ 
\mathcal{C}^{i}(V)\overset{\delta _{i}^{0}}{\longrightarrow }\mathcal{C}%
^{n-1+i}(V)\overset{\delta _{i}^{1}}{\longrightarrow }\mathcal{C}%
^{2(n-1)+i}(V)\longrightarrow \cdots \longrightarrow \mathcal{C}%
^{k(n-1)+i}(V)\overset{\delta _{i}^{k}}{\longrightarrow }\mathcal{C}%
^{(k+1)(n-1)+i}(V)\longrightarrow \cdots \\ 
\vdots \\ 
\mathcal{C}^{n-2}(V)\overset{\delta _{n-2}^{0}}{\longrightarrow }\mathcal{C}%
^{2n-3}(V)\overset{\delta _{n-2}^{1}}{\longrightarrow }\mathcal{C}%
^{3(n-1)-1}(V)\longrightarrow \cdots \longrightarrow \mathcal{C}%
^{(1+k)(n-1)-1}(V)\overset{\delta _{n-2}^{k}}{\longrightarrow }\mathcal{C}%
^{(k+2)(n-1)-1}(V)\longrightarrow \cdots \\ 
\end{array}%
\end{equation*}%
We can consider the associated cohomology.

\subsection{Second case: $n$ is odd}

Consider a $n$-ary multiplication \ associated to Gerstenhber's product with
an odd $n$. Then pre-Lie Identity applied to triples $(\varphi ,\mu ,\mu )$
with $\varphi \in \mathcal{C}^{ k}(V)$ is always fulfiled. To define a
cohomology for these algebras we have to restrict the space of cochains to
the subspace $\chi ^{k}(V)$ of $k$-linear applications $\varphi :V^{\otimes
k}\longrightarrow V$ subject to the following axioms 
\begin{equation*}
\left\{ 
\begin{array}{l}
(\varphi \bullet \mu )\bullet \mu =0, \\ 
(\mu \bullet \varphi )\bullet \mu =0, \\ 
\mu \bullet (\varphi \bullet \mu )=0.%
\end{array}%
\right.
\end{equation*}%
Pre-Lie identity applied to the triple $(\mu ,\varphi ,\mu )$ implies 
\begin{equation*}
(\mu \bullet \varphi )\bullet \mu -\mu \bullet (\varphi \bullet \mu )=(\mu
\bullet \mu )\bullet \varphi -\mu \bullet (\mu \bullet \varphi )
\end{equation*}%
so 
\begin{equation*}
(\mu \bullet \varphi )\bullet \mu =\mu \bullet (\varphi \bullet \mu )-\mu
\bullet (\mu \bullet \varphi ).
\end{equation*}%
If we moreover assume that $\varphi $ belongs to $\chi ^{k}(V)$ then $\mu
\bullet (\mu \bullet \varphi )=0.$

\begin{theorem}
Let 
\begin{equation*}
\partial ^{k}:\chi ^{k}(V)\longrightarrow \mathcal{C}^{k+n-1}(V)
\end{equation*}%
be the linear map defined by 
\begin{equation*}
\partial ^{k}\varphi =(-1)^{k-1}\mu \bullet \varphi -\varphi \bullet \mu .
\end{equation*}%
Then

1. The image of $\partial ^{k}$ is included in $\chi ^{k+n-1}(V)$.

2. We obtain the following identity 
\begin{equation*}
\partial ^{k+n-1}\circ \partial ^{k}=0.
\end{equation*}
\end{theorem}

\noindent \textit{Proof.} Let $\varphi $ be in $\chi ^{k}(V)$ and consider $%
\partial ^{k}\varphi .$ Then 
\begin{equation*}
(\partial ^{k}\varphi \bullet \mu )\bullet \mu =(-1)^{k-1}((\mu \bullet
\varphi )\bullet \mu )\bullet \mu -((\varphi \bullet \mu )\bullet \mu
)\bullet \mu =0,
\end{equation*}%
and 
\begin{equation*}
(\mu \bullet \partial ^{k}\varphi )\bullet \mu =(-1)^{k-1}(\mu \bullet (\mu
\bullet \varphi ))\bullet \mu -(\mu \bullet (\varphi \bullet \mu ))\bullet
\mu =0,
\end{equation*}%
finally 
\begin{equation*}
\mu \bullet (\partial ^{k}\varphi \bullet \mu )=(-1)^{k-1}\mu \bullet ((\mu
\bullet \varphi )\bullet \mu )-\mu \bullet ((\varphi \bullet \mu )\bullet
\mu )=0.
\end{equation*}%
Thus $\partial ^{k}\varphi \in \chi ^{k+n-1}(V)$. But 
\begin{equation*}
\begin{array}{ll}
(\partial ^{k+n-1}\circ \partial ^{k})\varphi & =\partial
^{k+n-1}((-1)^{k-1}\mu \bullet \varphi -\varphi \bullet \mu ) \\ 
& =\mu \bullet (\mu \bullet \varphi )+(-1)^{k}\mu \bullet (\varphi \bullet
\mu )+(-1)^{k}(\mu \bullet \varphi )\bullet \mu +(\varphi \bullet \mu
)\bullet \mu =0%
\end{array}%
\end{equation*}%
so 
\begin{equation*}
\partial ^{k+n-1}\circ \partial ^{k}=0
\end{equation*}%
which proves the result.

\medskip

\begin{corollary}
Considering $\delta _{i}^{j}=\partial ^{i+j(n-1)}$ we get the following
complexes: 
\begin{equation*}
\begin{array}{l}
\medskip \chi ^{0}(V)\overset{\delta _{0}^{0}}{\longrightarrow }\chi
^{n-1}(V)\overset{\delta _{0}^{1}}{\longrightarrow }\chi
^{2n-2}(V)\longrightarrow \cdots \longrightarrow \chi ^{k(n-1)}(V)\overset{%
\delta _{0}^{k}}{\longrightarrow }\chi ^{(k+1)(n-1)}(V)\longrightarrow \cdots
\\ 
\medskip \chi ^{1}(V)\overset{\delta _{1}^{0}}{\longrightarrow }\chi ^{n}(V)%
\overset{\delta _{1}^{1}}{\longrightarrow }\chi ^{2n-1}(V)\longrightarrow
\cdots \longrightarrow \chi ^{k(n-1)+1}(V)\overset{\delta _{1}^{k}}{%
\longrightarrow }\chi ^{(k+1)(n-1)+1}(V)\longrightarrow \cdots \\ 
\vdots \\ 
\chi ^{i}(V)\overset{\delta _{i}^{0}}{\longrightarrow }\chi ^{n-1+i}(V)%
\overset{\delta _{i}^{1}}{\longrightarrow }\chi ^{2n-1+i}(V)\longrightarrow
\cdots \longrightarrow \chi ^{k(n-1)+i}(V)\overset{\delta _{i}^{k}}{%
\longrightarrow }\chi ^{(k+1)(n-1)+i}(V)\longrightarrow \cdots \\ 
\vdots \\ 
\chi ^{n-2}(V)\overset{\delta _{n-2}^{0}}{\longrightarrow }\chi ^{2n-3}(V)%
\overset{\delta _{n-2}^{1}}{\longrightarrow }\chi
^{3(n-1)-1}(V)\longrightarrow \cdots \longrightarrow \chi ^{(1+k)(n-1)-1}(V)%
\overset{\delta _{n-2}^{k}}{\longrightarrow }\chi
^{(k+2)(n-1)-1}(V)\longrightarrow \cdots \\ 
\end{array}%
\end{equation*}
\end{corollary}

\subsection{Remark}

Let $(V,\mu )$ be an algebra of type $\bullet _{n,n}$. It is unital if there
exists $1\in V$ such that 
\begin{equation*}
\mu (1,1,\cdots ,X)=\mu (1,\cdots ,X,1)=\cdots =\mu (X,\cdots ,1)=X
\end{equation*}%
for any $X\in V$. Then for any $f\in End(V)$ we associate the bilinear
map $\varphi _{f}$ defined by: 
\begin{equation*}
\varphi _{f}(X,Y)=\partial ^{1}f(1,1,\cdots ,X,Y).
\end{equation*}%
Similary, for any bilinear application $\varphi $, we can associate the
trilinear application $\psi _{\varphi }$ given by : 
\begin{equation*}
\psi _{\varphi }(X,Y,Z)=\partial ^{2}\varphi (1,1,\cdots ,1,X,Y,Z)
\end{equation*}%
and if $\varphi $ belongs to $C^{k}(V)$ or $\chi ^{k}(V)$ we consider $\psi
_{\varphi }$ belonging to $\mathcal{C}^{k+1}(V)$ or $\chi ^{k+1}(V)$ given
by 
\begin{equation*}
\psi _{\varphi }(X_{1},\cdots ,X_{k+1})=\partial ^{k}\varphi (1,\cdots
,1,X_{1},\cdots ,X_{k+1}).
\end{equation*}%
Then we get the sequence 
\begin{equation*}
\mathcal{C}^{1}(V)\overset{\phi _{1}}{\longrightarrow }\mathcal{C}^{2}(V)%
\overset{\phi _{2}}{\longrightarrow }\mathcal{C}^{3}(V)\longrightarrow
\cdots \longrightarrow \mathcal{C}^{k}(V)\overset{\phi _{k}}{\longrightarrow 
}\mathcal{C}^{k+1}(V)\cdots
\end{equation*}%
where $\phi _{k}\varphi =\psi _{\varphi }.$ Computing $\phi _{k+1}\circ \phi
_{k}$ we get 
\begin{equation*}
\phi _{k+1}(\phi _{k}(\varphi ))=(\partial ^{k+1}(\partial ^{k}\varphi
))(1,\cdots ,1,X_{1},\cdots ,X_{k})=0.
\end{equation*}%
Thus the previous sequence is a complex and we get: 
\begin{equation*}
\begin{array}{llllll}
\medskip &  & \downarrow \phi _{n-2} &  & \downarrow \phi _{k(n-1)-1} &  \\ 
\medskip \mathcal{C}^{0}(V) & \overset{\delta _{0}^{0}}{\longrightarrow } & 
\mathcal{C}^{n-1}(V) & \overset{\delta _{0}^{1}}{\longrightarrow }\cdots
\longrightarrow & \mathcal{C}^{k(n-1)}(V) & \overset{\delta _{0}^{k}}{%
\longrightarrow }\cdots \\ 
\medskip \downarrow \phi _{0} &  & \downarrow \phi _{n-1} &  & \downarrow
\phi _{k(n-1)} &  \\ 
\mathcal{C}^{1}(V) & \overset{\delta _{1}^{0}}{\longrightarrow } & \mathcal{C%
}^{n}(V) & \overset{\delta _{1}^{1}}{\longrightarrow }\cdots \longrightarrow
& \mathcal{C}^{k(n-1)+1}(V) & \overset{\delta _{1}^{k}}{\longrightarrow }%
\cdots \\ 
\medskip \downarrow \phi _{1} &  & \downarrow \phi _{n-1+1} &  & \downarrow
\phi _{k(n-1)+1} &  \\ 
\vdots &  & \vdots &  & \vdots &  \\ 
\medskip \downarrow \phi _{i-1} &  & \downarrow \phi _{n-1+i-1} &  & 
\downarrow \phi _{k(n-1)+i-1} &  \\ 
\medskip \mathcal{C}^{i}(V) & \overset{\delta _{i}^{0}}{\longrightarrow } & 
\mathcal{C}^{n-1+i}(V) & \overset{\delta _{i}^{1}}{\longrightarrow }\cdots
\longrightarrow & \mathcal{C}^{k(n-1)+i}(V) & \overset{\delta _{i}^{k}}{%
\longrightarrow }\cdots \\ 
\downarrow \phi _{i} &  & \downarrow \phi _{n-1+i} &  & \downarrow \phi
_{k(n-1)+i} &  \\ 
\vdots &  & \vdots &  & \vdots &  \\ 
\downarrow \phi _{n-3} &  & \downarrow \phi _{n-1+n-3} &  & \downarrow \phi
_{(1+k)(n-1)-2} &  \\ 
\mathcal{C}^{n-2}(V) & \overset{\delta _{n-2}^{0}}{\longrightarrow } & 
\mathcal{C}^{2n-3}(V) & \overset{\delta _{n-2}^{1}}{\longrightarrow }\cdots
\longrightarrow & \mathcal{C}^{(1+k)(n-1)-1}(V) & \overset{\delta _{n-2}^{k}}%
{\longrightarrow }\cdots \\ 
\medskip \downarrow \phi _{n-2} &  & \downarrow \phi _{2n-3} &  & \downarrow
\phi _{(1+k)(n-1)-1} &  \\ 
\mathcal{C}^{n-1}(V) & \overset{\delta _{0}^{1}}{\longrightarrow } & 
\mathcal{C}^{2(n-1)}(V) & \overset{\delta _{0}^{2}}{\longrightarrow }\cdots
\longrightarrow & \mathcal{C}^{(k+1)(n-1)}(V) & \overset{\delta _{0}^{k+1}}{%
\longrightarrow }\cdots%
\end{array}%
\end{equation*}

\section{Deformations and cohomology}

Let $A=(V,\mu )$ be a $n$-ary partially associative algebra. By a
deformation of $(V,\mu )$ we mean a $\mathbb{K}[[t]]$-$n$-ary partially
associative algebra $A_{t}=(V_{t},\mu _{t})$ where $V_{t}=V\otimes \mathbb{K}%
[[t]]$ and $A_{t}/(tA_{t})\simeq A.$

We know that there exists always a cohomology theory which controls these
deformations. We recall the construction. Let $M$ be the variety of
structure constants and $\mathbb{K}[M]$ the affine coordinate ring of $M$.
We construct a resolution $(\Lambda (X)_{\ast },d)\rightarrow (\mathbb{K}%
[M],d=0)$ where $X$ is a graded vector space $X=\underset{i\geq 0}{\oplus }%
X_{i}$ and $\Lambda (X)_{\ast }$ a graded commutative algebra on $X$, the
differential $d$ satisfying $d(X_{i})\subset \Lambda (X_{i-1})$ and $%
H_{i}(\Lambda (X),d)=0$ for $i>1.$ If $L^{\ast }=Der(\Lambda (X)_{\ast
},\Lambda (X)_{\ast })$ and $\delta $ the differential on $L^{\ast }$
induced by $d,$ then $H^{\ast }(L,\delta )$ controls the deformations. But
this cohomology is too general to be useful for practical computations.

\noindent In our case, taking for any $\varphi \in \mathcal{C}^{k}(V),$ 
$$\delta^k\varphi=(-1)^{k-1}\mu\bullet_{2,k}\varphi-\varphi\bullet_{k,2}\mu,$$ we have a
complex 
$$(\overline{\chi }^{k}(V)=\chi^{k}(V)\oplus Ker\,\delta
^{k},\delta ^{k})$$
 such that $H^{2}(\overline{\chi }^{\ast },\overline{\chi }%
^{\ast })$ controls the deformations.

\section{Graded version of Gerstenhaber's products and associated $n$-ary
algebras}

\subsection{A relation of degree $7$}

In this section we claim that $n$ is a natural odd number. In this case we
already know that, for a cochain $\phi \in \mathcal{C}^{k}(V)$, the identity
of $n$-ary algebra $(V,\mu )$ 
\begin{equation*}
(\phi \bullet _{k,n}\mu )\bullet _{k+n-1,n}\mu =0
\end{equation*}%
is not always fulfiled (contrary to the even case) and that we must impose
that the cochains satisfy this identity to define a cohomology. As $\mu \bullet \mu =0$, this
identity is equivalent to: 
\begin{equation*}
\phi \circ \theta _{k}(\mu )=\phi \circ \underset{%
\begin{array}{c}
{\small 0\leq p\leq k-2} \\ 
{\small 0\leq q\leq k-2-p}%
\end{array}%
}{\sum }Id_{p}\otimes \mu \otimes Id_{q}\otimes \mu \otimes Id_{k-p-q-2}=0.
\end{equation*}%
In particular, as $\mu \bullet \mu =0$, we get 
\begin{equation*}
\mu \circ \theta _{k}(\mu )=0.
\end{equation*}%
This identity for $n=3$\ writes 
\begin{equation*}
\begin{array}{l}
\mu \circ (Id_{1}\otimes \mu \otimes Id_{1})\circ (\mu \otimes Id_{4})+\mu
\circ (Id_{2}\otimes \mu )\circ (\mu \otimes Id_{4})+\mu \circ
(Id_{2}\otimes \mu )\circ (Id_{1}\otimes \mu \otimes Id_{3}) \\ 
+\mu \circ (\mu \otimes Id_{2})\circ (Id_{3}\otimes \mu \otimes Id_{1})+\mu
\circ (\mu \otimes Id_{2})\circ (Id_{4}\otimes \mu )+\mu \circ
(Id_{1}\otimes \mu \otimes Id_{1})\circ (Id_{4}\otimes \mu )=0.%
\end{array}%
\end{equation*}%
so we get 
\begin{equation*}
\begin{array}{l}
(\mu \circ (Id_{1}\otimes \mu \otimes Id_{1})\circ (\mu \otimes Id_{4})+\mu
\circ (\mu \otimes Id_{2})\circ (Id_{3}\otimes \mu \otimes Id_{1})) \\ 
+(\mu \circ (Id_{2}\otimes \mu )\circ (\mu \otimes Id_{4})+\mu \circ (\mu
\otimes Id_{2})\circ (Id_{4}\otimes \mu )) \\ 
+(\mu \circ (Id_{2}\otimes \mu )\circ (Id_{1}\otimes \mu \otimes Id_{3})+\mu
\circ (Id_{1}\otimes \mu \otimes Id_{1})\circ (Id_{4}\otimes \mu ))=0.%
\end{array}%
\end{equation*}%
Similary the identity $\mu \bullet (\phi \bullet \mu )=0$ is equivalent to: 
\begin{equation*}
\begin{array}{l}
\mu \circ (Id_{1}\otimes \phi \otimes Id_{1})\circ (Id_{1}\otimes \mu
\otimes Id_{3}+Id_{2}\otimes \mu \otimes Id_{2}+Id_{3}\otimes \mu \otimes
Id_{1}) \\ 
+\mu \circ (\phi \otimes Id_{2})\circ (\mu \otimes Id_{4}+Id\otimes \mu
\otimes Id_{3}+Id_{2}\otimes \mu \otimes Id_{2}) \\ 
+\mu \circ (Id_{2}\otimes \phi )\circ (Id_{4}\otimes \mu +Id_{3}\otimes \mu
\otimes Id_{1}+Id_{2}\otimes \mu \otimes Id_{2})=0 \\ 
\end{array}%
\end{equation*}%
and for $\phi =\mu $ this identity is fulfiled. It writes 
\begin{equation*}
\begin{array}{l}
\mu \circ (\phi \otimes Id_{2})\circ (\mu \otimes Id_{4}-Id_{3}\otimes \mu
\otimes Id_{1})+\mu \circ (Id_{2}\otimes \phi )\circ (Id_{4}\otimes \mu
-Id_{1}\otimes \mu \otimes Id_{3})=0. \\ 
\end{array}%
\end{equation*}%
Then we get

\begin{proposition}
Let $(V,\mu )$ be a ternary algebra with multiplication of type $\bullet
_{3,3}.$ Then $\mu $ satisfies the following relations of degree $7$:

1%
${{}^\circ}$%
) 
\begin{equation*}
\begin{array}{l}
\mu\circ (Id_1\otimes \mu \otimes Id_1)\circ (\mu \otimes Id_4) +\mu\circ (
\mu \otimes Id_2)\circ (Id_3\otimes\mu \otimes Id_1)+ \mu\circ (Id_2\otimes
\mu )\circ (\mu \otimes Id_4) \\ 
+\mu\circ ( \mu \otimes Id_2)\circ (Id_4\otimes\mu ) +\mu\circ (Id_2\otimes
\mu )\circ (Id_1\otimes \mu \otimes Id_3) \mu\circ (Id_1\otimes \mu \otimes
Id_1)\circ (Id_4 \otimes \mu)=0.%
\end{array}%
\end{equation*}
\medskip

2%
${{}^\circ}$%
) 
\begin{equation*}
\begin{array}{l}
\mu\circ ( \mu \otimes Id_2)\circ (\mu \otimes Id_4 -Id_3\otimes\mu \otimes
Id_1) +\mu\circ (Id_2\otimes \mu )\circ (Id_4\otimes\mu - Id_1\otimes\mu
\otimes Id_3)=0. \\ 
\end{array}%
\end{equation*}
\end{proposition}

The interpretation of the first relation show the necessity to distinguish
the order of multiplications. A classical meethod consists in grading the
initail space so we get: 
\begin{equation*}
(I)\ \ \ \left\{ 
\begin{array}{l}
\mu \circ (Id_{1}\otimes \mu \otimes Id_{1})\circ (\mu \otimes Id_{4})=-\mu
\circ (\mu \otimes Id_{2})\circ (Id_{3}\otimes \mu \otimes Id_{1}) \\ 
\mu \circ (Id_{2}\otimes \mu )\circ (\mu \otimes Id_{4})=-\mu \circ (\mu
\otimes Id_{2})\circ (Id_{4}\otimes \mu ) \\ 
\mu \circ (Id_{2}\otimes \mu )\circ (Id_{1}\otimes \mu \otimes Id_{3})=-\mu
\circ (Id_{1}\otimes \mu \otimes Id_{1})\circ (Id_{4}\otimes \mu ).%
\end{array}%
\right.
\end{equation*}%
We will shortely develop this approach in the following.

\subsection{Graded identities}

For any two maps $f\in \medskip \mathcal{C}^{k}(V)$ and $g\in \medskip 
\mathcal{C}^{l}(V)$ we consider 
\begin{equation*}
f\bullet _{i}g(X_{1},\cdots ,X_{k+l-1})=f(X_{1},\cdots
,X_{i-1},g(X_{i},\cdots ,X_{i+l-1}),\cdots ,X_{k+l-1})
\end{equation*}%
so 
\begin{equation*}
f\bullet _{k,l}g=\overset{k}{\underset{i=1}{\sum }}(-1)^{(i-1)(l-1)}f\bullet
_{i}g
\end{equation*}

We will now work in a $\mathbb{Z}$-graded vector space $V=\oplus _{n\in 
\mathbb{Z}}V_{n}$. We define the suspension (resp. desuspension) of $V$ by $%
\uparrow V$ (resp. $\downarrow V$), i.e. the graded $\mathbb{Z}$-graded
vector space $\uparrow V=\oplus _{n\in \mathbb{Z}}(\uparrow V)_{n}$ $\ $%
(resp. $\downarrow V=\oplus _{n\in \mathbb{Z}}(\downarrow V)_{n}$) with $%
(\uparrow V)_{n}=V_{n+1}$(resp. $(\downarrow V)_{n}=V_{n-1}$). So the
corresponding degree $+1$ map $\uparrow :V\longrightarrow \ \uparrow V$
sends $v\in V$ into its suspended copy $\uparrow $ $v\in \ \uparrow V,$
assigns to $V$ the graded vector space $\uparrow V$ and satisfies 
\begin{equation*}
\uparrow \circ \downarrow =\downarrow \circ \uparrow =Id.
\end{equation*}%
More generally we have 
\begin{equation*}
\uparrow ^{\otimes l}\circ \downarrow ^{\otimes l}=\downarrow ^{\otimes
l}\circ \uparrow ^{\otimes l}=(-1)^{l(l-1)/2}Id
\end{equation*}

\noindent Suppose that the algebra $V$ is graded. If $f:V^{\otimes
k}\longrightarrow V$ has a degree $|f|$, then if 
\begin{equation*}
\phi (f)=\uparrow \circ f\circ \downarrow ^{\otimes k}
\end{equation*}%
we get 
\begin{equation*}
\phi (f)\bullet _{i}\phi (g)=(-1)^{(|g|+k-1)(l-i)+|g|(i-1)}\phi (f\bullet
_{i}g)
\end{equation*}%
for graded $f\in C^{k}(A)$ and $g\in C^{l}(A)$. Let $\mu \in C^{n}(A)$ be an
application of degree $n-2$. We get 
\begin{equation*}
\phi (\mu )\bullet _{i}\phi (\mu )=(-1)^{(n-2+n-1)(n-i)+(n-2)(i-1)}\phi (\mu
\bullet _{i}\mu )=(-1)^{i(n+1)}\phi (\mu \bullet _{i}\mu ).
\end{equation*}%
Thus 
\begin{equation*}
\begin{array}{ll}
\medskip \phi (\mu \bullet _{n,n}\mu )=\phi
(\sum_{i=1}^{n}(-1)^{(i-1)(n-1)}\mu \bullet _{i}\mu ) & 
=\sum_{i=1}^{n}(-1)^{(i-1)(n-1)}\phi \mu \bullet _{i}\mu ) \\ 
\medskip & =\sum_{i=1}^{n-1}(-1)^{(n-1)}\phi (\mu )\bullet _{i}\phi (\mu )
\\ 
\medskip & =(-1)^{(n-1)}\sum_{i=1}^{n-1}\phi (\mu )\bullet _{i}\phi (\mu ).%
\end{array}%
\end{equation*}%
For example for $n=3$, the graded identity $\mu \bullet _{3,3}\mu $ writes 
\begin{equation*}
\phi (\mu \bullet _{3,3}\mu )=\sum \phi (\mu )\bullet _{i}\phi (\mu )
\end{equation*}%
et pour $n=2$ 
\begin{equation*}
\phi (\mu \bullet _{2,2}\mu )=-\phi (\mu )\bullet _{1}\phi (\mu )-\phi (\mu
)\bullet _{2}\phi (\mu )
\end{equation*}%
All these identites are sign constant. In particular:

\begin{theorem}
Let $V=\oplus _{n\in \mathbb{Z}}V_{n}$ be a $\mathbb{Z}$-graded vector
space. A graded application $\mu $ with degree $n-2$ is a Gerstenhaber
multiplication of type $\bullet _{n,n}$ if and only if 
\begin{equation*}
\sum \phi (\mu )\bullet _{i}\phi (\mu )=0.
\end{equation*}
\end{theorem}

\subsection{Composition relations}

In the nongraded case we have: 
\begin{equation*}
\left\{ 
\begin{array}{l}
(\mu \bullet _{j}\mu )\bullet _{i}\mu =(\mu \bullet _{i}\mu )\bullet
_{j+n-1}\mu \ \ {\mbox{\rm if }}\ \ i+1\leq 2n-1 \\ 
(\mu \bullet _{j}\mu )\bullet _{i}\mu =(\mu \bullet _{i+n-1}\mu )\bullet
_{j}\mu \ \ {\mbox{\rm if}}\ \ 1\leq j\leq i-n\ \ {\mbox{\rm and}}\ \ i\geq
n+1%
\end{array}%
\right.
\end{equation*}
If $\mu $ is graded with degree $|\mu |$, the commutative rules come from
Koszul signs conventions 
\begin{equation*}
\left\{ 
\begin{array}{l}
(\mu \bullet _{j}\mu )\bullet _{i}\mu =(-1)^{|\mu ||\mu |}(\mu \bullet
_{i}\mu )\bullet _{j+n-1}\mu \ \ {\mbox{\rm if }}\ \ i+1\leq 2n-1, \\ 
(\mu \bullet _{j}\mu )\bullet _{i}\mu =(-1)^{|\mu ||\mu |}(\mu \bullet
_{i+n-1}\mu )\bullet _{j}\mu \ \ {\mbox{\rm if}}\ \ 1\leq j\leq i-n\ \ {%
\mbox{\rm et}}\ \ i\geq n+1.%
\end{array}%
\right.
\end{equation*}

\noindent \textit{Fundamental examples.}

\medskip

i) For $n=2$, $\mu $ is of degree $0$ and we obtain the relations of the non
graded case.

\medskip

ii) For $n=3$, $\mu $ is of degree $1$ and we get the relations

\begin{equation*}
\left\{ 
\begin{array}{l}
(\mu \bullet _{2}\mu )\bullet _{1}\mu =-(\mu \bullet _{1}\mu )\bullet
_{4}\mu, \\ 
(\mu \bullet _{3}\mu )\bullet _{1}\mu =-(\mu \bullet _{1}\mu )\bullet
_{5}\mu, \\ 
(\mu \bullet _{3}\mu )\bullet _{2}\mu =-(\mu \bullet _{2}\mu )\bullet
_{5}\mu.%
\end{array}%
\right.
\end{equation*}%
This gives us the relations claimed in (\textit{I})

\subsection{On the cohomology in the graded case for $n=3$}

Let $\uparrow A$ be the suspension of the graded space $A$. Consider $\mu $
as an application of degree $1$ 
\begin{equation*}
\mu : (\uparrow  A)^{\otimes 3}\longrightarrow (\uparrow  A)
\end{equation*}

\begin{proposition}
If $\mu $ is a $3$-ary multiplication of degree $1$, any $\varphi \in 
\mathcal{C}^{n}(\uparrow  A)$ satisfies 
\begin{equation*}
(\varphi \bullet \mu )\bullet \mu =0.
\end{equation*}
\end{proposition}

\noindent \textit{Consequence.}

\bigskip

Let $\delta :C^{n}( \uparrow A)\longrightarrow C^{n+2}(\uparrow A)$ be the $1$
degree operation defined by 
\begin{equation*}
\delta \varphi =\mu \bullet \varphi -(-1)^{|\varphi |}\varphi \bullet \mu
\end{equation*}%
where $|\varphi |$ is the degree of $\varphi .$

\begin{lemm}
(Graded pre-Lie identity)

Let $\varphi _{1}$be in $C^{n}(\uparrow A)$,$\varphi _{2}$ in $%
C^{m}(\uparrow A)$ and $\varphi _{3}$ in $C^{p}(\uparrow A).$ They satisfy 
\begin{equation*}
(\varphi _{1}\bullet \varphi _{2})\bullet \varphi _{3}-\varphi _{1}\bullet
(\varphi _{2}\bullet \varphi _{3})=(-1)^{(m-1)(p-1)}(-1)^{|\varphi
_{2}||\varphi _{3}|}((\varphi _{1}\bullet \varphi _{3})\bullet \varphi
_{2}-\varphi _{1}\bullet (\varphi _{3}\bullet \varphi _{2}))
\end{equation*}
\end{lemm}

We deduce 
\begin{equation*}
(\mu \bullet \mu )\bullet \varphi -\mu \bullet (\mu \bullet \varphi
)=(-1)^{|\varphi |}((\mu \bullet \varphi )\bullet \mu -\mu \bullet (\varphi
\bullet \mu ))
\end{equation*}%
and 
\begin{equation*}
\delta (\delta \varphi )=0
\end{equation*}

\begin{proposition}
The operator $\delta :C^{n}(\uparrow A)\longrightarrow C^{n+2}(\uparrow A)$
defined by 
\begin{equation*}
\delta \varphi =\mu \bullet \varphi -(-1)^{|\varphi |}\varphi \bullet \mu
\end{equation*}%
gives the complex 
\begin{equation*}
C^{0}(\uparrow A)\longrightarrow C^{3}(\uparrow A)\longrightarrow \cdots
\end{equation*}
\end{proposition}

We denote $H^{\ast }(\uparrow A,\delta {\mu })$ the associated cohomology.

\medskip

\noindent \textbf{Remark.} In [M,R] we give an explanation in operadic terms
as the underlying quadratic operad are not Koszul with the usual definition
of the cohomology $H^{\ast }(\uparrow A,\delta {\mu }).$

\section{The free algebra $L(V,$ $\bullet _{3,3})$}

Let $V$ be a $\mathbb{K}$-vector space. For even $k,$ the free algebras $%
L(V,_{{}}\bullet _{k,k})$ have been described in \cite{Gne}. But the odd
case behaves in a completely different way, as we have already seen it for
the cohomology. We are going to describe in detail the case $k=3$ that is, the
case of a $3$-ary algebra $V$ with multiplication $\bullet _{3,3}$.

As we have a $3$ order product, the free algebra is graded as follows%
\begin{equation*}
L(V,_{{}}\bullet _{3,3})=\oplus _{p\geq 1}L^{2p+1}(V)
\end{equation*}%
with%
\begin{equation*}
L^{1}(V)=V,\quad L^{3}(V)=V^{\otimes ^{3}}.
\end{equation*}
Let's describe the further terms.%
\begin{equation*}
L^{5}(V)=((V^{\otimes ^{3}}\otimes V^{\otimes ^{2}})\oplus (V\otimes
V^{\otimes ^{3}}\otimes V)\oplus (V^{\otimes ^{2}}\otimes V^{\otimes
^{3}}))/R_{5}
\end{equation*}%
where $R_{5}$ is the sub-space of $(V^{\otimes ^{3}}\otimes V^{\otimes
^{2}})\oplus (V\otimes V^{\otimes ^{3}}\otimes V)\oplus (V^{\otimes
^{2}}\otimes V^{\otimes ^{3}})$ of relations spanned with vectors which write%
\begin{equation*}
(v_{1}\otimes v_{2}\otimes v_{3})\otimes v_{4}\otimes v_{5}+v_{1}\otimes
(v_{2}\otimes v_{3}\otimes v_{4})\otimes v_{5}+v_{1}\otimes v_{2}\otimes
(v_{3}\otimes v_{4}\otimes v_{5}).
\end{equation*}%
If $V$ is $n$-dimensional $\dim L^{5}(V)=2n^{5}.$

To describe the other components we denote by $D(k,3)$, for any positif odd $%
k,$ the set of triples $(a,b,c)$ satisfying 
\begin{equation*}
\left\{ 
\begin{array}{l}
a,b,c\in \text{odd positif integers}, \\ 
a+b+c=k.%
\end{array}%
\right.
\end{equation*}%
We also need to use a simplified notation for vectors replacing a term $%
v_{i_1}\otimes v_{i_2}\otimes \cdots \otimes v_{i_p}$ by $i_1 i_2 \cdots i_p$
For example the vector $(v_{1}\otimes v_{2}\otimes v_{3})\otimes
v_{4}\otimes v_{5}$ writes $(1\cdot 2\cdot 3)\cdot 4\cdot 5.$ We now consider%
\begin{equation*}
L^{7}(V)=(\oplus _{(a,b,c)\in D(7,3)}L^{a}(V)\otimes L^{b}(V)\otimes
L^{c}(V))/R_{7}
\end{equation*}%
where $R_{7}$ is the sub-space of $\oplus _{(a,b,c)\in
D(7,3)}L^{a}(V)\otimes L^{b}(V)\otimes L^{c}(V)$ spanned with the vectors%
\begin{equation*}
\left\{ 
\begin{array}{c}
((1 \cdot 2 \cdot 3)\cdot 4\cdot 5) \cdot 6 \cdot 7+(1\cdot 2\cdot 3)\cdot
(4\cdot 5\cdot 6)\cdot 7 +(1\cdot 2\cdot 3)\cdot 4\cdot (5\cdot 6\cdot 7),
\\ 
(1\cdot (2\cdot 3\cdot 4)\cdot 5)\cdot 6\cdot 7+1\cdot ((2\cdot 3\cdot
4)\cdot 5\cdot 6)\cdot 7 +1\cdot (2\cdot 3\cdot 4)\cdot (5\cdot 6\cdot 7),
\\ 
(1\cdot 2\cdot (3\cdot 4\cdot 5))\cdot 6\cdot 7+1\cdot (2\cdot (3\cdot
4\cdot 5)\cdot 6)\cdot 7 +1\cdot 2\cdot ((3\cdot 4\cdot 5)\cdot 6\cdot 7)),
\\ 
(1\cdot 2\cdot 3)\cdot (4\cdot 5\cdot 6)\cdot 7+1\cdot ((2\cdot 3\cdot
(4\cdot 5\cdot 6))\cdot 7 +1\cdot 2\cdot (3\cdot (4\cdot 5\cdot 6)\cdot 7),
\\ 
((1\cdot 2\cdot 3)\cdot 4\cdot (5\cdot 6\cdot 7)+1\cdot (2\cdot 3\cdot
4)\cdot (5\cdot 6\cdot 7) +1\cdot 2\cdot (3\cdot 4\cdot (5\cdot 6\cdot 7)).%
\end{array}%
\right.
\end{equation*}%
As $\dim (\oplus _{(a,b,c)\in D(7,3)}L^{a}(V)\otimes L^{b}(V)\otimes
L^{c}(V) $ )$=3(\dim L^{5}(V)\times (\dim V)^{2}+3(\dim L^{3}(V))^{2}\times
\dim V)=9n^{7},$ we deduce that $\dim L^{7}(V)=4n^{7}.$ \ To describe the
general case we need a more efficient coding of vectors. An element of $%
L^{2p+1}(V)$ writes $1\cdot (2(\cdots (\cdots ))\cdot 2p+1)$ with $p-1$
brackets. As each bracket has to contain $3$ elements, we can code an
element of $L^{2p+1}(V)$ by the position of the left brackets. For example $%
1\cdot 2\cdot (3\cdot 4\cdot (5\cdot 6\cdot 7))$\ corresponds to $g_{3,5}$
as the left brackets are at the element $3$ and $5$. We also suppose that we
point the brackets from the left to the right, that is, $g_{j_{1},\cdots ,
j_{p-1}}$ belongs to $L^{2p+1}$ with $j_{1}\leq j_{2}\leq ...\leq
j_{p-2}\leq j_{p-1}$. Thus the above elements of $L^{7}$ correspond to $%
g_{1,1}, \, g_{1,2},\, g_{1,3},\, g_{1,4},\, g_{1,5},\, g_{2,2},\,
g_{2,3},\, g_{2,4},\, g_{2,5},\, g_{3,3},\, g_{3,4}$ and $g_{3,5}$ and $%
R_{7} $ are the relations spanned by the vectors 
\begin{equation*}
\left\{ 
\begin{array}{l}
g_{1,1}+g_{1,4}+g_{1,5}, \\ 
g_{1,2}+g_{2,2}+g_{2,5}, \\ 
g_{1,3}+g_{2,3}+g_{3,3}, \\ 
g_{1,4}+g_{2,4}+g_{3,4}, \\ 
g_{1,5}+g_{2,5}+g_{3,5}, \\ 
g_{1,1}+g_{1,2}+g_{1,3}, \\ 
g_{2,2}+g_{2,3}+g_{2,4}, \\ 
g_{3,3}+g_{3,4}+g_{3,5}.%
\end{array}
\right.
\end{equation*}

This coding allows to define the sub-space of relations in any degree by an
inductive way. An element of $L^{2p+1}$ has the coding $g_{j_{1},\cdots ,
j_{p-1}}$ with $1\leq j_{1}\leq 3,$ $j_{1}\leq j_{2}\leq 5, \cdots ,
j_{p-2}\leq j_{p-1}\leq 2p-1.$ Suppose that we have the relations $%
R_{2p-1}.$ These relations concern the vectors coded by $g_{j_{1},\cdots ,\,
j_{p-2}}.\ $\ The relations of $R_{2p+1}$ are obtained by relations of $%
R_{2p-1}$ using two rules: Suppose we have a relation in $R_{2p-1}$, that is
implying vectors $g_{j_{1},\cdots ,\, j_{p-2}}.$ We have to explain how we
get from such a vector $g_{j_{1},\cdots ,\, j_{p-2}}$ a vector $g_{\,
l_{1},\cdots ,\, l_{p-1}}$ involved in $R_{2p-1}$

\begin{itemize}
\item Add the index $i$ in front of each $(p-2)$-uple of the vectors $%
g_{j_{1},\cdots ,j_{p-2}}$ involved in the relation, with $i$ successively
equal to $1,2~$\ and $3$. We replace each index $j_{l}$ by $j_{l}+(i-1).$

For example: \ $g_{1,4}$ becomes successively $g_{1,1,4},\, g_{2,2,5}$ and $%
g_{3,3,6}.$

\item For $i$ successively equal to $1,2,\cdots ,2p-1,$ add the index $i$ in
front of any $(p-2)$-uple of vector $g_{j_{1},\cdots ,\,j_{p-2}}$ involved
in the relation; if the index $j_{1}$ is less than or equal to $i$, conserve 
$j_{1}$, otherwise replace $j_{1}$ par $j_{1}+2.$ And do the same for all
further indices.\ Rearrange the indices to get $1\leq j_{1}\leq 3,$ $%
j_{1}\leq j_{2}\leq 5,\cdots ,\, j_{p-2}\leq j_{p-1}\leq 2p-1.$
\end{itemize}

Thus any relations in $R_{2p-1}$ gives $(2p-1)+3=2p+2$ \ relations in $%
R_{2p+1}.$ We have then constructed the generating relations of $R_{2p+1}.$

\bigskip

\noindent \textbf{Example : relations of }$R_{9}.$ Each of the $8$ relations
of $R_{7}$ gives $10$ relations.\ For example $g_{1,1}+g_{1,2}+g_{1,3}$
gives 
\begin{equation*}
\left\{ 
\begin{array}{l}
g_{1,1,1}+g_{1,1,2}+g_{1,1,3}, \\ 
g_{2,2,2}+g_{2,2,3}+g_{2,2,4} , \\ 
g_{3,3,3}+g_{3,3,4}+g_{3,3,5}, \\ 
g_{1,1,1}+g_{1,1,4}+g_{1,1,5}, \\ 
g_{1,1,2}+g_{1,2,2}+g_{1,2,5}, \\ 
g_{1,1,3}+g_{1,2,3}+g_{1,3,3}, \\ 
g_{1,1,4}+g_{1,2,4}+g_{1,3,4}, \\ 
g_{1,1,5}+g_{1,2,5}+g_{1,3,5}, \\ 
g_{1,1,6}+g_{1,2,6}+g_{1,3,6}, \\ 
g_{1,1,7}+g_{1,2,7}+g_{1,3,7}.%
\end{array}
\right.
\end{equation*}
We then get $80$ relations. We can solve this system directely or using
computer.\ We solved this system using Mathematica and found $\dim
R_{9}=20n^{9}$ (the rank of the system ist $20).\ $Thus $\dim
L^{9}(V)=5n^{9}.$

\bigskip

\noindent \textbf{Remark.} Contrary to the previous case there exists some
trivial homogeneous products. Any element of $L^{9}$ is considered as a
product of $3$ elements, i.e $u\in L^9, $ $u=u_a \otimes u_b \otimes u_c$
with $(a,b,c) \in D(9,3)$ and $u_a \in L^a,\, u_b \in L^b, \, u_c \in L^c,$
or more simply, $u$ is of type $(a,b,c) \in D(9,3).$ All the elements having
a factor in $L^{7}$ (that is of type $(3,3,1),(3,1,3),(1,3,3)$) are zero.
Also all the homogeneous products of type $(5,3,1),\, (3,5,1),\, (1,5,3),\,
(1,3,5)$ whose elements in $5$ elements are of type $(113)31,\, 3(311)1,\,
1(113)3)$ and $13(311)$ are zero. By elements of type $(113)31$ we mean
elements which can be written $(x_{i_1} \otimes x_{i_2} \otimes x_{i_3})
\otimes x_{i_4} \otimes x_5$ with $x_{i_1},x_{i_2},x_{i_5}\in V$ and $%
x_{i_3},x_{i_4} \in L^3(V)$. At least the elements of type $11(11(113)),\,
11((311)11),\, 1(11(113))1,\, (11(113))11, \, ((311)11)11$ and $(333)$ are
zero. The rule defining these elements is the following: let us consider an
element as a product of $\ 3$ elements $(a,b,c).$ Thus the elements
containing

\begin{itemize}
\item $3$ products of $L^{3}$,

\item $2$ products of $L^{3}$ in the same bracket (for example $(3,1,3)$),

\item $2$ products of $L^{3}$ consecutif but in different brackets (for
example $1(113)3$),

\item only $1$ product of $L^{3}$ but neighboring to $2$ brackets (for
example $11((311)11)$
\end{itemize}
are all zero. Moreover remark that a basis of $L^{9^{{}}}$ is given by the
vectors coded by 
\begin{equation*}
g_{3,4,4},\, g_{3,4,6},\, g_{1,2,4},\, g_{1,2,2},\, g_{1,1,7}.
\end{equation*}

\begin{definition}
Let $V$ be a vector space. The free algebra of type $\bullet _{3,3}$ on $V$
is the $3$-ary algebra $L(V,_{{}}\bullet _{3,3})=\oplus _{p\geq
1}L^{2p+1}(V) $ with $L^{1}(V)=V,\quad L^{3}(V)=V^{\otimes ^{3}}$and%
\begin{equation*}
L^{2p+1}(V)=(\oplus _{(a,b,c)\in D(2p+1,3)}L^{a}(V)\otimes L^{b}(V)\otimes
L^{c}(V))/R_{2p+1}
\end{equation*}%
where $R_{2p+1}$ is the sub-space of $\oplus _{(a,b,c)\in
D(2p+1,3)}L^{a}(V)\otimes L^{b}(V)\otimes L^{c}(V)$ spanned by the vectors $%
g_{j_{1},\cdots ,j_{p-1}\text{ }}$with $1\leq j_{1}\leq 3,$ $j_{1}\leq
j_{2}\leq 5,\cdots \, , j_{p-2}\leq j_{p-1}\leq 2p-1$ and satisfying the
relations defined by the above rules.
\end{definition}

It is clear that $L(V,_{{}}\bullet _{3,3})=\oplus _{p\geq 1}L^{2p+1}(V)$ is
of type $\bullet _{3,3}.\ $If $a_{1},a_{2},a_{3}$ are three homogeneous
elements, $a_{i}\in L^{2p_{i}+1}$, the product is defined by the class of $%
a_{1}\otimes a_{2}\otimes a_{3}.$ This algebra satisfies the following
property:

\begin{proposition}
Let $\mathcal{A}$ be a $3$-ary algebra of type $\bullet _{3,3}$ and $V$ a
vector space. Then any linear map $f:V\rightarrow \mathcal{A}$ can be
factorized in a unique morphism of $3$-ary algebras%
\begin{equation*}
F:L(V,\bullet _{3,3})\rightarrow \mathcal{A}\text{.}
\end{equation*}
\end{proposition}
\textit{Proof.} If $(\mathcal{A}_{1},\mu _{1})$ and $(\mathcal{A}%
_{2},\mu _{2})$ are $3$-ary algebras, a linear map $g:\mathcal{A}%
_{1}\rightarrow \mathcal{A}_{2}$ is a morphism of algebras if%
\begin{equation*}
\mu _{2}(g(X),g(Y),g(Z))=g(\mu _{1}(X,Y,Z)),
\end{equation*}%
for any $X,Y,Z\in \mathcal{A}_{1}.$ Let $f:V\rightarrow \mathcal{A}$ be a
linear map. Consider the linear map $F:L(V,\bullet _{3,3})\rightarrow 
\mathcal{A}$ defined on homogeneous components of $L^{2p+1}(V)$ by 
\begin{equation*}
F(g_{j_{1},\cdots ,\, j_{p-1}}\otimes (v_{p}\otimes \cdots \otimes
v_{2p+1}))=g_{j_{1},\cdots ,\, j_{p-1}}\otimes (f(v_{p})\otimes \cdots
\otimes f(v_{2p+1})),
\end{equation*}%
where\ $g_{j_{1},\cdots ,\, j_{p-1}}\otimes (v_{p}\otimes \cdots \otimes v_{2p+1})$
corresponds to the vector $(v_{1}\otimes v_{2}\otimes \cdots \otimes
(v_{j_{1}}\otimes \cdots \otimes (v_{j_{2}}\otimes \cdots \otimes
(v_{j_{p-1}}\otimes v_{j_{p}}\otimes v_{j_{p+1}})\cdots ).$ We obtain the
expected morphism of algebras.

It remains to give a basis of the free algebra. We have already computed the
dimensions of the first homogeneous components. Let us complete these
results by describing a basis. For this we use a graphic representation by a
planary trees with $3$-branching noods (three entries and one exit as the
multiplication is $3$-ary).\ $\ $\ We decorate each leave with a basis
vector of $V$ to obtain a free familly of elements of the free algebra.
Suppose $V$ is $n$-dimensional$.$ Then

\noindent $\bullet $ $\dim L^{3}(V)=n^{3}.$ A basis is associated to the
tree 
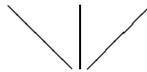
\begin{equation*}
\xymatrix{ && \\ &\ar@{-}[ul]\ar@{-}[u]\ar@{-}[ur]& \\ \\ }
\end{equation*}%
$\bullet $ $\dim L^{5}(V)=2n^{5}.$ A basis is associated to the trees%
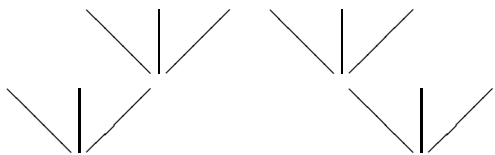
\begin{equation*}
\xymatrix{ &&& \\ &&\ar@{-}[ul]\ar@{-}[u]\ar@{-}[ur]& \\ &
\ar@{-}[ul]\ar@{-}[u]\ar@{-}[ur]&& \\ \\ }\ \xymatrix{ &&&& \\
&\ar@{-}[ul]\ar@{-}[u]\ar@{-}[ur]&&& \\ &&
\ar@{-}[ul]\ar@{-}[u]\ar@{-}[ur]&& \\ \\ }
\end{equation*}%
$\bullet $ $\dim L^{7}(V)=4n^{7}.$ A basis corresponds to the trees\noindent 
{\tiny 
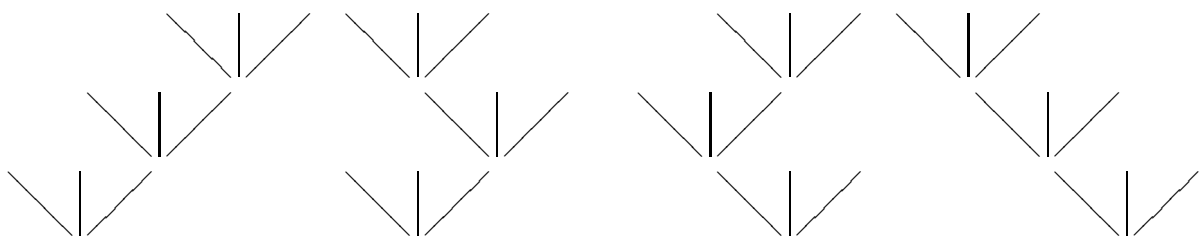
\begin{equation*}
\xymatrix{ &&&& \\ &&&\ar@{-}[ul]\ar@{-}[u]\ar@{-}[ur]& \\ &&
\ar@{-}[ul]\ar@{-}[u]\ar@{-}[ur]&& \\ & \ar@{-}[ul]\ar@{-}[u]\ar@{-}[ur]&&&
\\ }\ \xymatrix{ &&&& \\ &\ar@{-}[ul]\ar@{-}[u]\ar@{-}[ur]&&& \\ &&
\ar@{-}[ul]\ar@{-}[u]\ar@{-}[ur]&& \\ & \ar@{-}[ul]\ar@{-}[u]\ar@{-}[ur]&&&
\\ }\!\!\!\!\!\!\!\!\!\!\!\!\!\!\!\!\!\!\!\!\!\!\!\!\!\!\!\!\!\!\!\!\!\!\!\!%
\!\xymatrix{ &&&& \\ &&&\ar@{-}[ul]\ar@{-}[u]\ar@{-}[ur]& \\ &&
\ar@{-}[ul]\ar@{-}[u]\ar@{-}[ur]&& \\ &&&
\ar@{-}[ul]\ar@{-}[u]\ar@{-}[ur]&\\ }\ \xymatrix{ &&&& \\
&\ar@{-}[ul]\ar@{-}[u]\ar@{-}[ur]&&& \\ &&
\ar@{-}[ul]\ar@{-}[u]\ar@{-}[ur]&& \\ &&& \ar@{-}[ul]\ar@{-}[u]\ar@{-}[ur]&
\\ }
\end{equation*}%
} $\bullet $ $\dim L^{9}(V)=5n^{9}.$ A basis is given by\noindent 
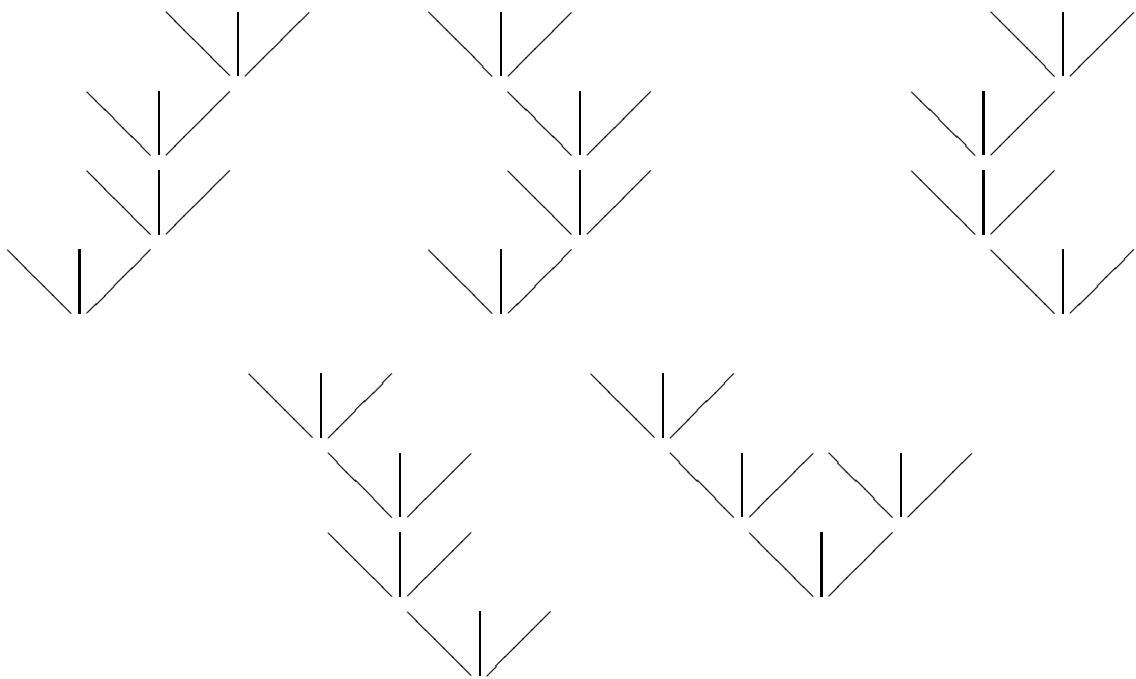
\begin{equation*}
\xymatrix{ &&&&& \\ &&&\ar@{-}[ul]\ar@{-}[u]\ar@{-}[ur]&& \\ &&
\ar@{-}[ul]\ar@{-}[u]\ar@{-}[ur]&&& \\ &&
\ar@{-}[ul]\ar@{-}[u]\ar@{-}[ur]&&& \\ & \ar@{-}[ul]\ar@{-}[u]\ar@{-}[ur]&&&
}\ \xymatrix{ &&&&& \\ &\ar@{-}[ul]\ar@{-}[u]\ar@{-}[ur]&&& \\ &&
\ar@{-}[ul]\ar@{-}[u]\ar@{-}[ur]&&& \\ &&
\ar@{-}[ul]\ar@{-}[u]\ar@{-}[ur]&&& \\ & \ar@{-}[ul]\ar@{-}[u]\ar@{-}[ur]&&&
}\!\!\!\!\!\!\!\!\!\!\!\!\!\!\!\!\!\!\!\!\xymatrix{ &&&&& \\
&&&&\ar@{-}[ul]\ar@{-}[u]\ar@{-}[ur]& \\ &&&
\ar@{-}[ul]\ar@{-}[u]\ar@{-}[ur]&& \\ &&&
\ar@{-}[ul]\ar@{-}[u]\ar@{-}[ur]&&\\ &&&&
\ar@{-}[ul]\ar@{-}[u]\ar@{-}[ur]&\\ }
\end{equation*}%
\begin{equation*}
\xymatrix{ &&&&& \\ &&\ar@{-}[ul]\ar@{-}[u]\ar@{-}[ur]&&& \\ &&&
\ar@{-}[ul]\ar@{-}[u]\ar@{-}[ur]&& \\ &&&
\ar@{-}[ul]\ar@{-}[u]\ar@{-}[ur]&&\\ &&&&
\ar@{-}[ul]\ar@{-}[u]\ar@{-}[ur]&\\ }\ \xymatrix{ &&&&& \\
&\ar@{-}[ul]\ar@{-}[u]\ar@{-}[ur]&&&& \\ &&
\ar@{-}[ul]\ar@{-}[u]\ar@{-}[ur]&&\ar@{-}[ul]\ar@{-}[u]\ar@{-}[ur]& \\ &&&
\ar@{-}[ul]\ar@{-}[u]\ar@{-}[ur]&& \\ }
\end{equation*}%
$\bullet $ $\dim L^{11}(V)=6n^{11}.$ A basis is given by\noindent 
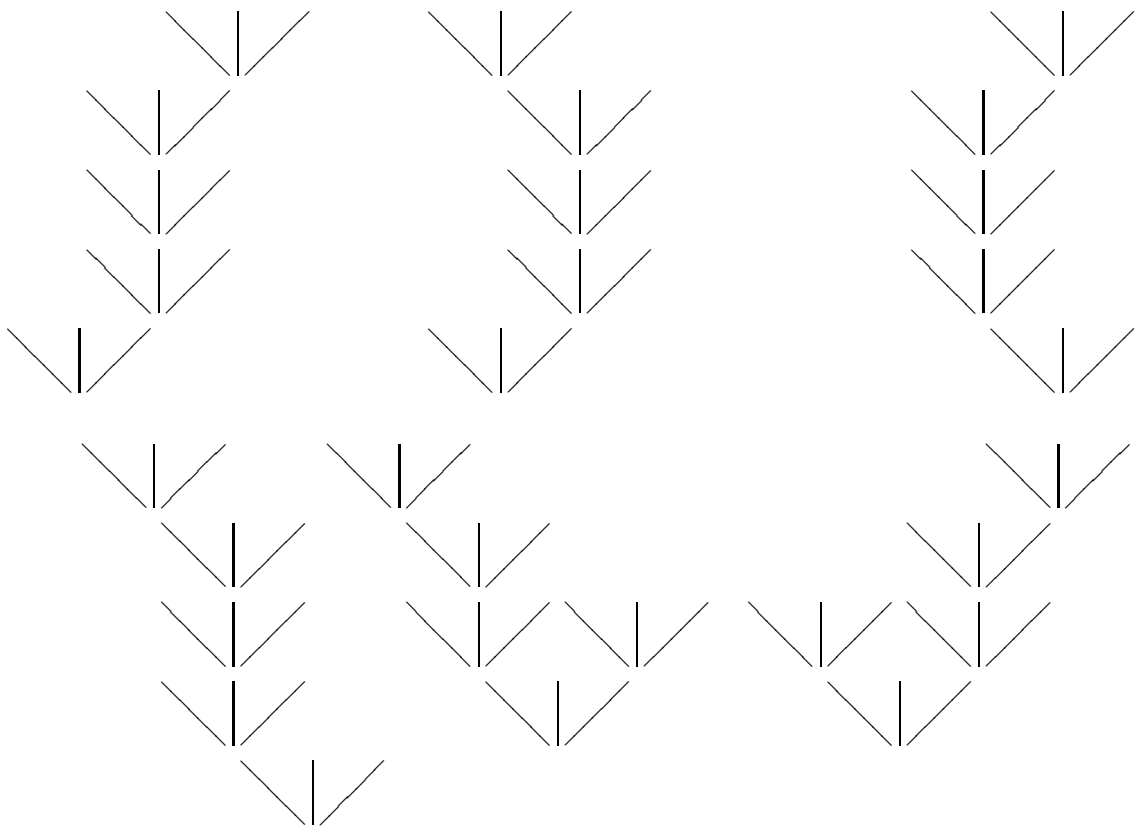
\begin{equation*}
\xymatrix{ &&&&& \\ &&&\ar@{-}[ul]\ar@{-}[u]\ar@{-}[ur]&& \\ &&
\ar@{-}[ul]\ar@{-}[u]\ar@{-}[ur]&&& \\ &&
\ar@{-}[ul]\ar@{-}[u]\ar@{-}[ur]&&& \\ &&
\ar@{-}[ul]\ar@{-}[u]\ar@{-}[ur]&&& \\ & \ar@{-}[ul]\ar@{-}[u]\ar@{-}[ur]&&&
}\ \xymatrix{ &&&&& \\ &\ar@{-}[ul]\ar@{-}[u]\ar@{-}[ur]&&& \\ &&
\ar@{-}[ul]\ar@{-}[u]\ar@{-}[ur]&&& \\ &&
\ar@{-}[ul]\ar@{-}[u]\ar@{-}[ur]&&& \\ &&
\ar@{-}[ul]\ar@{-}[u]\ar@{-}[ur]&&& \\ & \ar@{-}[ul]\ar@{-}[u]\ar@{-}[ur]&&&
}\!\!\!\!\!\!\!\!\!\!\!\!\!\!\!\!\!\!\!\!\xymatrix{ &&&&& \\
&&&&\ar@{-}[ul]\ar@{-}[u]\ar@{-}[ur]& \\ &&&
\ar@{-}[ul]\ar@{-}[u]\ar@{-}[ur]&& \\ &&& \ar@{-}[ul]\ar@{-}[u]\ar@{-}[ur]&&
\\ &&& \ar@{-}[ul]\ar@{-}[u]\ar@{-}[ur]&&\\ &&&&
\ar@{-}[ul]\ar@{-}[u]\ar@{-}[ur]&\\ }
\end{equation*}%
\begin{equation*}
\xymatrix{ &&&&& \\ &&\ar@{-}[ul]\ar@{-}[u]\ar@{-}[ur]&&& \\ &&&
\ar@{-}[ul]\ar@{-}[u]\ar@{-}[ur]&& \\ &&& \ar@{-}[ul]\ar@{-}[u]\ar@{-}[ur]&&
\\ &&& \ar@{-}[ul]\ar@{-}[u]\ar@{-}[ur]&&\\ &&&&
\ar@{-}[ul]\ar@{-}[u]\ar@{-}[ur]&\\ }\!\!\!\!\!\!\!\!\!\!\!\!\!\!\!\!\!\!\!\!%
\xymatrix{ &&&&& \\ &\ar@{-}[ul]\ar@{-}[u]\ar@{-}[ur]&&& \\ &&
\ar@{-}[ul]\ar@{-}[u]\ar@{-}[ur]&&&\\ &&
\ar@{-}[ul]\ar@{-}[u]\ar@{-}[ur]&&\ar@{-}[ul]\ar@{-}[u]\ar@{-}[ur]& \\ &&&
\ar@{-}[ul]\ar@{-}[u]\ar@{-}[ur]&& \\ }\ \xymatrix{ &&&&& \\
&&&&\ar@{-}[ul]\ar@{-}[u]\ar@{-}[ur]& \\
&&&\ar@{-}[ul]\ar@{-}[u]\ar@{-}[ur]&& \\ &\ar@{-}[ul]\ar@{-}[u]\ar@{-}[ur]
&&\ar@{-}[ul]\ar@{-}[u]\ar@{-}[ur]&& \\ &&
\ar@{-}[ul]\ar@{-}[u]\ar@{-}[ur]&&& }\ 
\end{equation*}
$\bullet $ $\dim L^{13}(V)=7n^{13}.$ A basis is given by\noindent 
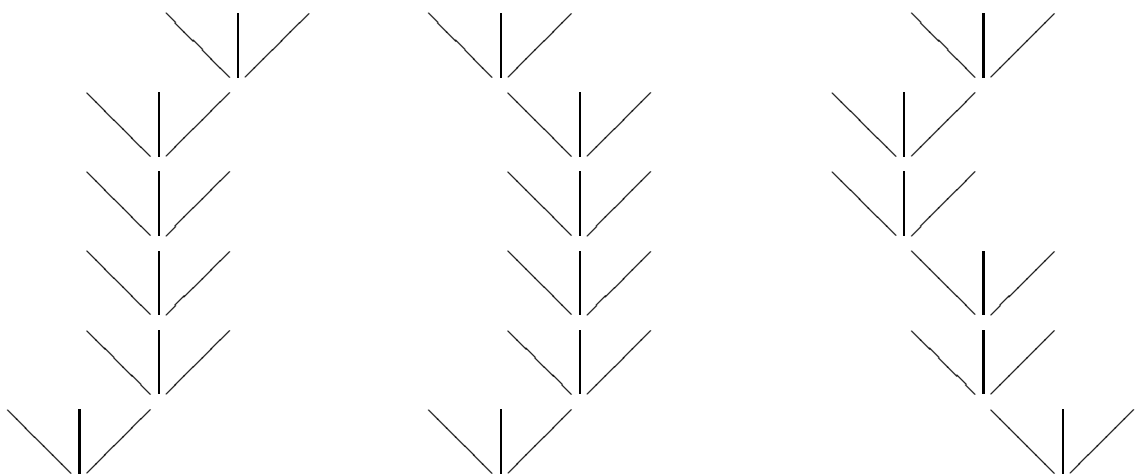
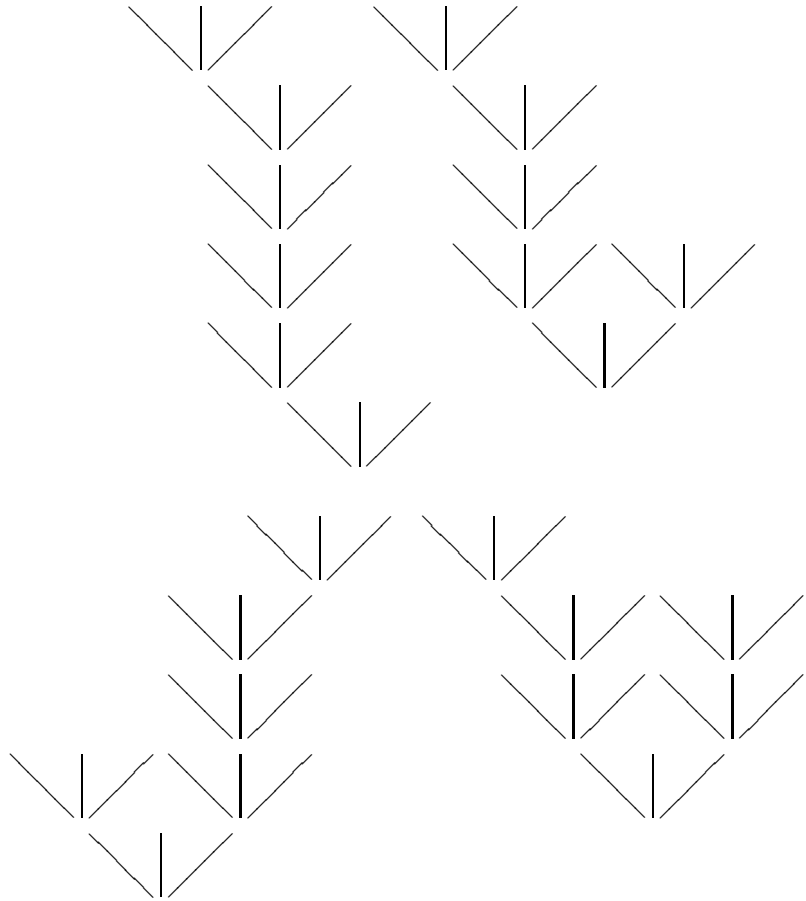
\begin{equation*}
\xymatrix{ &&&&& \\ &&&\ar@{-}[ul]\ar@{-}[u]\ar@{-}[ur]&& \\ &&
\ar@{-}[ul]\ar@{-}[u]\ar@{-}[ur]&&& \\ &&
\ar@{-}[ul]\ar@{-}[u]\ar@{-}[ur]&&& \\ &&
\ar@{-}[ul]\ar@{-}[u]\ar@{-}[ur]&&& \\ &&
\ar@{-}[ul]\ar@{-}[u]\ar@{-}[ur]&&& \\ & \ar@{-}[ul]\ar@{-}[u]\ar@{-}[ur]&&&
}\ \xymatrix{ &&&&& \\ &\ar@{-}[ul]\ar@{-}[u]\ar@{-}[ur]&&& \\ &&
\ar@{-}[ul]\ar@{-}[u]\ar@{-}[ur]&&& \\ &&
\ar@{-}[ul]\ar@{-}[u]\ar@{-}[ur]&&& \\ &&
\ar@{-}[ul]\ar@{-}[u]\ar@{-}[ur]&&& \\ &&
\ar@{-}[ul]\ar@{-}[u]\ar@{-}[ur]&&& \\ & \ar@{-}[ul]\ar@{-}[u]\ar@{-}[ur]&&&
}\!\!\!\!\!\!\!\!\!\!\!\!\!\!\!\!\!\!\!\!\xymatrix{ &&&&& \\
&&&\ar@{-}[ul]\ar@{-}[u]\ar@{-}[ur]& \\ &&
\ar@{-}[ul]\ar@{-}[u]\ar@{-}[ur]&& \\ && \ar@{-}[ul]\ar@{-}[u]\ar@{-}[ur]&&
\\ &&& \ar@{-}[ul]\ar@{-}[u]\ar@{-}[ur]&& \\ &&&
\ar@{-}[ul]\ar@{-}[u]\ar@{-}[ur]&&\\ &&&&
\ar@{-}[ul]\ar@{-}[u]\ar@{-}[ur]&\\ }
\end{equation*}%
\begin{equation*}
\xymatrix{ &&&&& \\ &&\ar@{-}[ul]\ar@{-}[u]\ar@{-}[ur]&&& \\ &&&
\ar@{-}[ul]\ar@{-}[u]\ar@{-}[ur]&& \\ &&& \ar@{-}[ul]\ar@{-}[u]\ar@{-}[ur]&&
\\ &&& \ar@{-}[ul]\ar@{-}[u]\ar@{-}[ur]&& \\ &&&
\ar@{-}[ul]\ar@{-}[u]\ar@{-}[ur]&&\\ &&&&
\ar@{-}[ul]\ar@{-}[u]\ar@{-}[ur]&\\ }\!\!\!\!\!\!\!\!\!\!\!\!\!\!\!\!\!\!\!\!%
\xymatrix{ &&&&& \\ &\ar@{-}[ul]\ar@{-}[u]\ar@{-}[ur]&&&& \\
&&\ar@{-}[ul]\ar@{-}[u]\ar@{-}[ur]&&&\\ &&
\ar@{-}[ul]\ar@{-}[u]\ar@{-}[ur]&&&\\ &&
\ar@{-}[ul]\ar@{-}[u]\ar@{-}[ur]&&\ar@{-}[ul]\ar@{-}[u]\ar@{-}[ur]& \\ &&&
\ar@{-}[ul]\ar@{-}[u]\ar@{-}[ur]&& \\ }\ 
\end{equation*}
\begin{equation*}
\xymatrix{ &&&&& \\ &&&& \ar@{-}[ul]\ar@{-}[u]\ar@{-}[ur]& \\ &&&
\ar@{-}[ul]\ar@{-}[u]\ar@{-}[ur]&& \\ &&& \ar@{-}[ul]\ar@{-}[u]\ar@{-}[ur]&&
\\ & \ar@{-}[ul]\ar@{-}[u]\ar@{-}[ur]&&\ar@{-}[ul]\ar@{-}[u]\ar@{-}[ur]&& \\
&& \ar@{-}[ul]\ar@{-}[u]\ar@{-}[ur]&&& \\ } \xymatrix{ &&&&& \\
&\ar@{-}[ul]\ar@{-}[u]\ar@{-}[ur]&&&& \\ &&
\ar@{-}[ul]\ar@{-}[u]\ar@{-}[ur]&&\ar@{-}[ul]\ar@{-}[u]\ar@{-}[ur]& \\ &&
\ar@{-}[ul]\ar@{-}[u]\ar@{-}[ur]&&\ar@{-}[ul]\ar@{-}[u]\ar@{-}[ur]& \\ &&&
\ar@{-}[ul]\ar@{-}[u]\ar@{-}[ur]&& \\ }
\end{equation*}
The choice of the basis is non canonical. But we choose them for symmetry
reasons. The rules providing the relations of the sub-space $R_{2p+1}$ are
easy to implement in order to solve the corresponding linear system. This
gives the dimensions of the spaces $L^{2p+1}(V)$ (in fact we find the
dimensions of the modules of the associated operad). We have illustrated
this approach in small dimensions above. Let us notice that we can however
present basic vectors for the relations associtated to the elements of 
\begin{equation*}
L^{2p-1}\otimes L^{1}\otimes L^{1}\oplus L^{1}\otimes L^{2p-1}\otimes
L^{1}\oplus L^{1}\otimes L^{1}\otimes L^{2p-1}.
\end{equation*}%
These elements correspond to the trees \noindent 
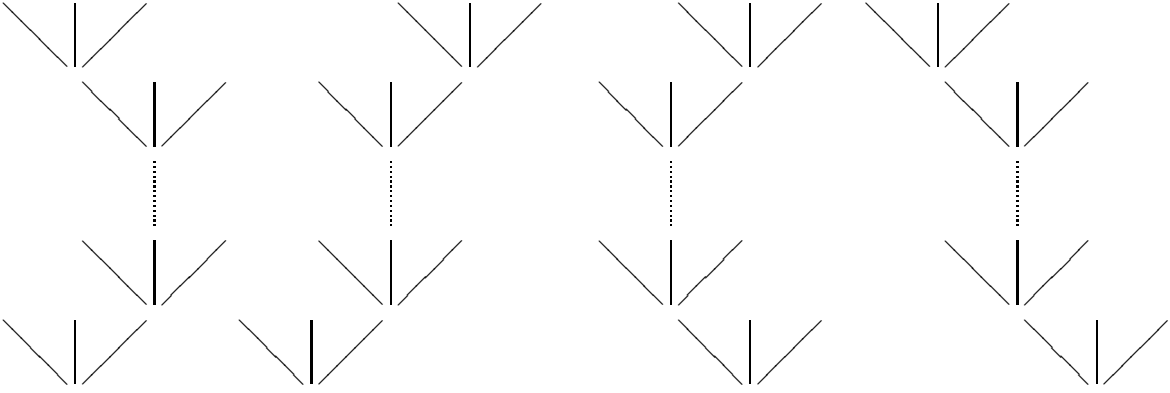
\begin{equation*}
\xymatrix{ &&&&& \\ &\ar@{-}[ul]\ar@{-}[u]\ar@{-}[ur]&&&& \\ &&
\ar@{-}[ul]\ar@{-}[u]\ar@{-}[ur]&&& \\ && \ar@{.}[u] &&& \\ &&
\ar@{-}[ul]\ar@{-}[u]\ar@{-}[ur]&&& \\ & \ar@{-}[ul]\ar@{-}[u]\ar@{-}[ur]&&&
}\!\!\!\!\!\!\!\!\!\!\!\!\!\!\!\!\!\!\!\!\!\!\!\!\!\!\!\!\!\!\!\!\!\!\!\!\!%
\!\!\!\xymatrix{ &&&&& \\ &&&\ar@{-}[ul]\ar@{-}[u]\ar@{-}[ur]&& \\ &&
\ar@{-}[ul]\ar@{-}[u]\ar@{-}[ur]&&& \\ && \ar@{.}[u] &&& \\ &&
\ar@{-}[ul]\ar@{-}[u]\ar@{-}[ur]&&& \\ & \ar@{-}[ul]\ar@{-}[u]\ar@{-}[ur]&&&
}\!\!\!\!\!\!\!\!\!\!\!\!\!\!\!\!\!\!\!\!\!\!\!\!\!\!\!\!\!\!\!\!\!\!\!\!\!%
\!\!\!\!\!\!\!\!\!\!\!\xymatrix{ &&&&& \\
&&&&\ar@{-}[ul]\ar@{-}[u]\ar@{-}[ur]& \\ &&&
\ar@{-}[ul]\ar@{-}[u]\ar@{-}[ur]&& \\ &&& \ar@{.}[u]&& \\ &&&
\ar@{-}[ul]\ar@{-}[u]\ar@{-}[ur]&&\\ &&&&
\ar@{-}[ul]\ar@{-}[u]\ar@{-}[ur]&\\ }\!\!\!\!\!\!\!\!\!\!\!\!\!\!\!%
\xymatrix{ &&&&& \\ && \ar@{-}[ul]\ar@{-}[u]\ar@{-}[ur]&&& \\ &&&
\ar@{-}[ul]\ar@{-}[u]\ar@{-}[ur] && \\ &&& \ar@{.}[u]&& \\ &&&
\ar@{-}[ul]\ar@{-}[u]\ar@{-}[ur]&&\\ &&&&
\ar@{-}[ul]\ar@{-}[u]\ar@{-}[ur]&\\ }
\end{equation*}
The other are of the form

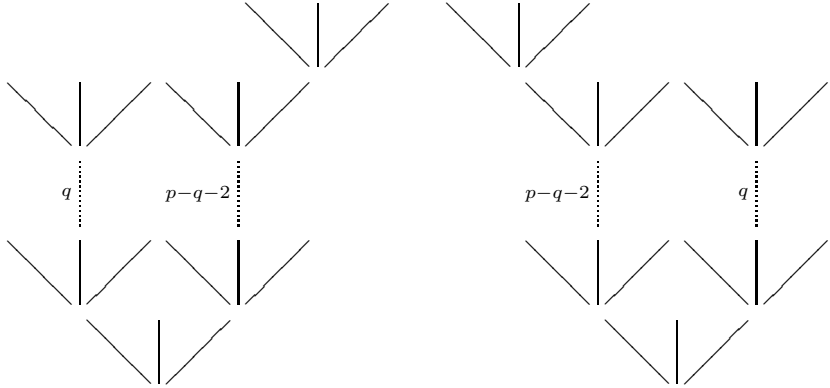
\begin{equation*}
\xymatrix{ &&&&&& \\ &&&&\ar@{-}[ul]\ar@{-}[u]\ar@{-}[ur]&& \\ &  \ar@{-}[ul]\ar@{-}[u]\ar@{-}[ur] &&
\ar@{-}[ul]\ar@{-}[u]\ar@{-}[ur]&&& \\ & \ar@{.}^{q}[u]  && \ar@{.}^{p-q-2}[u] &&& \\
 &\ar@{-}[ul]\ar@{-}[u]\ar@{-}[ur] && \ar@{-}[ul]\ar@{-}[u]\ar@{-}[ur]&&& \\ 
&& \ar@{-}[ul]\ar@{-}[u]\ar@{-}[ur] &&&
} \!\!\!\!\!\!\!\!\!\!\!\!\!\!\!\!\!\!\!\!\!\!\!\!\!\!\!\!\!\! \xymatrix{
&&&&&& \\ && \ar@{-}[ul]\ar@{-}[u]\ar@{-}[ur] &&&& \\ &&&
\ar@{-}[ul]\ar@{-}[u]\ar@{-}[ur] &&\ar@{-}[ul]\ar@{-}[u]\ar@{-}[ur] & \\ &&& \ar@{.}^{p-q-2}[u] &&\ar@{.}^{q}[u]& \\ &&&
\ar@{-}[ul]\ar@{-}[u]\ar@{-}[ur]&&\ar@{-}[ul]\ar@{-}[u]\ar@{-}[ur]& \\ &&&&
\ar@{-}[ul]\ar@{-}[u]\ar@{-}[ur]&&\\ }
\end{equation*}
where $q=1,\cdots,[\frac{p-2}{2}]$ and $[,]$ indicates the integer part of a rational number. If $p$ is even, the last two trees are
related. If $p$ is odd, these trees are independent. We deduce: 
\begin{theorem}
For any $p$ we have
$$dim L^{2p+1}(V)=(p+1)n^{2p+1}$$
where $n=dim V$.
\end{theorem}

\noindent \textbf{Remark.} Recall that for any vector space $V$, the
associated tensor algebra $T(V)$ is the unique solution, up to isomorphism,
of the universal problem which determine from a linear application $%
f:M\rightarrow A $ in an associative algebra $A$, a morphism of associative
algebra $T(V)\rightarrow A$. The construction of this algebra comes from the
isomorphisms%
\begin{equation*}
\Phi _{n,m}:T^{\otimes n}(V)\otimes T^{\otimes m}(V)\rightarrow T^{\otimes
(n+m)}(V)
\end{equation*}%
defined by%
\begin{equation*}
\Phi _{n,m}((x_{1}\otimes x_{2}\cdots \otimes x_{n})\otimes (y_{1}\otimes
y_{2}\cdots \otimes y_{m}))=x_{1}\otimes x_{2}\cdots x_{n}\otimes
y_{1}\otimes y_{2}\cdots \otimes y_{m}.
\end{equation*}%
In fact the multiplication $\mu $ of $T(V)$ is given by

\begin{equation*}
\mu ((x_{1}\otimes x_{2}\cdots \otimes x_{n})\otimes (y_{1}\otimes
y_{2}\cdots \otimes y_{m}))=\Phi _{n,m}((x_{1}\otimes x_{2}\cdots \otimes
x_{n})\otimes (y_{1}\otimes y_{2}\cdots \otimes y_{m}))
\end{equation*}%
and the associativity of the multiplication follows from 
\begin{equation*}
\Phi _{n+m,p}\bullet (\Phi _{n,m}\otimes Id_{p})=\Phi _{n+m,p}\bullet
(Id_{n}\otimes \Phi _{m,p}).
\end{equation*}%
We can define an other isomorphism non longer adapted to the associative
structure but adapted to the $n$-ary structure. For this we consider the familly of vectorial isomorphisms

\begin{equation*}
\Psi _{n,m,p}:T^{\otimes n}(V)\otimes T^{\otimes m}(V)\otimes T^{\otimes
p}(V)\rightarrow T^{\otimes n+m+p}(V)
\end{equation*}%
satisfying

\begin{equation*}
\left\{ 
\begin{array}{cc}
\Psi _{n,m+p+q,r}\bullet (Id_{n}\otimes \Psi _{m,p,q}\otimes Id_{r}) & 
=-2\Psi _{n,m+p+q,r}\bullet (Id_{n+m}\otimes \Psi _{p,q,r}) \\ 
& =-2\Psi _{n,m+p+q,r}\bullet (\Psi _{n,m,p}\bullet Id_{q+r}).
\end{array}%
\right.
\end{equation*}

\section{Extension of the notion of coassociative algebras for $n$-ary
algebras  }
 For $n=2$ we have that $\ 2$-ary partially associative algebras are
just associative algebras and we can define coassociative coalgebras with
the wellknown relations between these two structures. In fact, the dual
space of a coassociative algebra can be provided with a structure of
associative algebra, the dual space of a finite dimensional associative
algebra can be provided with a structure of coassociative coalgebra
structure and also, if $\ (A,\mu )$ is an associative algebra and $(M,\Delta
)$ a coassociative coalgebra, the space $Hom(M,A)$ can be provided with an
associative algebra structure. All these notions can be extended to $n$-ary algebras.

\noindent A $n$-ary partially associative algebra has a product $\mu $ satisfying
Equation (\ref{equa}) written in the following form%
$$
\sum_{p=0}^{n-1}  (-1)^{p(n-1)}\mu \circ (Id_{p}\otimes
\mu \otimes Id_{n-1-p})=0.
$$
Then we get the definition of partially coassociative $n$-ary
coalgebra.

\begin{definition}
A $n$-ary comultiplication on a $\mathbb{K}$-vector space $M$ is a map 
$$\Delta :M\rightarrow M^{\otimes ^{n}}.$$
A $n$-ary partially coassocative coalgebra  is a $\mathbb{K}$-vector space $%
M$ provided with a $n$-ary comultiplication $\Delta $ satisfying 
\begin{equation*}
{\overset{n-1}{\underset{p=0}{\sum }}}(-1)^{p(n-1)}(Id_{p}%
\otimes \Delta \otimes Id_{n-1-p})\circ \Delta =0.
\end{equation*}%
A $n$-ary totally coassocative coalgebra \ is a $\mathbb{K}$-vector space $M$
provided with a $n$-ary comultiplication $\Delta $ satisfying 
\begin{equation*}
(Id_{p}\otimes \Delta \otimes Id_{n-1-p})\circ \Delta =(Id_{q}\otimes
\Delta \otimes Id_{n-1-q})\circ \Delta ,
\end{equation*}%
for any $p,q\in \{0,\cdots ,n-1\}.$
\end{definition}
If $(\mathcal{A},\mu)$ is a $n$-ary algebra and $(M,\Delta)$ $n$-ary coalgebra
we denote by
$$
\begin{array}{ll}
A(\mu)=\underset{}{\overset{n-1}{\underset{p=0}{\sum }}}(-1)^{p(n-1)}\mu \circ (Id_{p}%
\otimes \mu \otimes Id_{n-1-p}), \\

\tilde{A}(\Delta)=
\underset{}{\overset{n-1}{\underset{p=0}{\sum }}}(-1)^{p(n-1)}(Id_{p}%
\otimes \Delta \otimes Id_{n-1-p})\circ \Delta.
\end{array}
$$
For any natural number $n$ and any $\mathbb{K}$-vector spaces $E$ and $F$,
we denote by 
\begin{equation*}
\lambda _{n}:Hom(E,F)^{\otimes n}\longrightarrow Hom(E^{\otimes
n},F^{\otimes n})
\end{equation*}%
the natural embedding 
\begin{equation*}
\lambda _{n}(f_{1}\otimes ...\otimes f_{n})(x_{1}\otimes ...\otimes
x_{n})=f_{1}(x_{1})\otimes ...\otimes f_{n}(x_{n}).
\end{equation*}

\begin{proposition}
The dual space of a $n$-ary partially coassociative coalgebra is provided
with a structure of $n$-ary partially associative algebra.
\end{proposition}
\textit{Proof}. Let $(M,\Delta )$ be a $n$-ary partially coassociative
coalgebra. We consider the multiplication on the dual vector space $M^{\ast }
$ of $M$ defined by : 
\begin{equation*}
\mu =\Delta ^{\ast }\circ \lambda _{n}.
\end{equation*}%
It provides $M^{\ast }$ with a $n$-ary partially assocative algebra
structure. In fact we have 
\begin{eqnarray} \label{equ}
\mu (f_{1}\otimes f_{2}\otimes \cdots \otimes f_{n})=\mu _{\mathbb{K}}\circ
\lambda _{n}(f_{1}\otimes f_{2}\otimes \cdots \otimes f_{n})\circ \Delta 
\end{eqnarray}%
for all $f_{1},\cdots ,f_{n}\in M^{\ast }$ where $\mu _{\mathbb{K}}$ is the
multiplication of $\mathbb{K}$. Equation (\ref{equ}) becomes : 
\begin{equation*}
\begin{array}{l}
\smallskip
\mu \circ (Id_{p}\otimes \mu \otimes Id_{n-1-p})(f_{1}\otimes f_{2}\otimes
\cdots \otimes f_{2n-1}) \\
\smallskip
=\mu _{\mathbb{K}}\circ (\lambda
_{n}(f_{1}\otimes \cdots \otimes f_{p}\otimes \mu (f_{p+1}\otimes \cdots
\otimes f_{p+n})\otimes f_{p+n+1}\otimes \cdots \otimes f_{2n-1}))\circ
\Delta  \\ 
\smallskip    \\ 
=\mu _{\mathbb{K}}\circ \lambda _{n}(f_{1}\otimes \cdots \otimes
f_{p}\otimes (\mu _{\mathbb{K}}\circ \lambda _{n}(f_{p+1}\otimes \cdots
\otimes f_{p+n})\circ \Delta )\otimes f_{p+n+1}\otimes \cdots \otimes
f_{2n-1})\circ \Delta  \\ 
\smallskip   \\ 
=\mu _{\mathbb{K}}\circ (Id_{p}\otimes \mu _{\mathbb{K}}\otimes
Id_{n-1-p})\circ \lambda _{2n-1}(f_{1}\otimes \cdots \otimes f_{2n-1})\circ
(Id_{p}\otimes \Delta \otimes Id_{n-1-p})\circ \Delta .%
\end{array}%
\end{equation*}%
Using associativity and commutativity of the multiplication in $\mathbb{K}$,
we obtain%
\begin{equation*}
\forall p,q\in \{0,\cdots ,n-1\},\quad \mu _{\mathbb{K}}\circ (Id_{p}\otimes \mu _{%
\mathbb{K}}\otimes Id_{n-1-p})=\mu _{\mathbb{K}}\circ (Id_{q}\otimes \mu _{%
\mathbb{K}}\otimes Id_{n-1-q}),
\end{equation*}%
so 
\begin{equation*}
\begin{array}{l}
\smallskip
\sum\limits_{p=0}^{n-1} (-1)^{p(n-1)}
\mu \circ (Id_{p}\otimes
\mu \otimes Id_{n-1-p}) \\
\smallskip
= \mu _{\mathbb{K}}\circ (\mu _{\mathbb{K}}\otimes
Id_{n-1})\circ \lambda _{2n-1}(f_{1}\otimes \cdots \otimes f_{2n-1})\circ 
\sum\limits_{p=0}^{n-1}(-1)^{p(n-1)}(Id_{p}\otimes \Delta
\otimes Id_{n-1-p})\circ \Delta  = 0 
\end{array}%
\end{equation*}%
and $(M^{\ast },\mu )$ is a $n$-ary partially partially associative algebra. 
$\Box $

\medskip

\begin{proposition}
The dual vector space of a finite dimensional $n$-ary partially associative
algebra has a $n$-ary partially associative coalgebra structure. 
\end{proposition}
\textit{Proof.} Let $\mathcal{A}$ be a finite dimensional $n$-ary partially
associative algebra and let $\{{e_{i}},{i=1,...,n}\}$ be a basis of $%
\mathcal{A}$. If $\{f_{i}\}$ is the dual basis then $\{f_{i_{1}}\otimes
\cdots \otimes f_{i_{n}}\}$ is a basis of $(\mathcal{A}^{\ast })^{\otimes n}$%
. The coproduct $\Delta $ on $\mathcal{A}^{\ast }$ is defined by 
\begin{equation*}
\Delta (f)=\sum_{i_{1},\cdots ,i_{n}}f(\mu (e_{i_{1}}\otimes \cdots \otimes
e_{i_{n}}))f_{i_{1}}\otimes \cdots \otimes f_{i_{n}}.
\end{equation*}%
In particular 
\begin{equation*}
\Delta (f_{k})=\sum_{i_{1},\cdots ,i_{n}}C_{i_{1}\cdots
,i_{n}}^{k}f_{i_{1}}\otimes \cdots \otimes f_{i_{n}}
\end{equation*}%
where $C_{i_{1}\cdots ,i_{n}}^{k}$ are the structure constants of $\mu $
related to the basis $\{{e_{i}}\}$. Then $\Delta $ is the comultiplication
of a $n$-ary partially coassocitive coalgebra. $\Box $

\medskip

Now we study the convolution product. Let us recall that if $(\mathcal{A},\mu )$ is associative $\mathbb{K}$%
-algebra and $(M,\Delta )$ a coassociative $\mathbb{K}$-coalgebra then the
convolution product 
\begin{equation*}
f\star g=\mu \circ \lambda _{2}(f\otimes g)\circ \Delta 
\end{equation*}%
provides $Hom(M,\mathcal{A})$ with an associative algebra structure. This
result can be extended to the $n$-ary partially associative algebras and
partially coassociative coalgebras. 

\begin{proposition}
Let $(\mathcal{A},\mu )$ be a $n$-ary partially associative algebra and $%
(M,\Delta )$ a $n$-ary totally coalgebra. Then the algebra $(Hom(M,\mathcal{A%
}),\star )$ is a $n$-ary partially associative algebra where $\star $ is the
convolution product : 
\begin{equation*}
f_{1}\star f_{2}\star \cdots \star f_{n}=\mu \circ \lambda _{n}(f_{1}\otimes
f_{2}\otimes \cdots \otimes f_{n})\circ \Delta .
\end{equation*}
\end{proposition}

\noindent \textit{Proof.} Let us compute the convolution product of  functions of $%
Hom(M,\mathcal{A})$. We have
\begin{equation*}
\begin{array}{l}
f_{1}\star \cdots \star f_{i-1}\star (f_{i}\star f_{i+1}\star \cdots \star
f_{i+n-1})\star f_{i+n}\star \cdots \star f_{2n-1} \smallskip \\
 =\mu \ \circ \lambda
_{n}(f_{1}\otimes f_{2}\otimes \cdots \otimes f_{i-1}\otimes (f_{i} \star \cdots \star f_{i+n-1})\otimes f_{i+n}\otimes \cdots \otimes f_{n})\,\circ \,\Delta  \\ 
  \vspace{0.1cm}\\ 
 =\mu \,\circ \,\lambda _{n}(f_1 \otimes \cdots f_{i-1}\otimes  (\mu \circ \lambda _{n}(f_{i}\otimes \cdots \otimes 
f_{i+n-1})\circ \Delta )\otimes f_{i+n} \otimes f_{2n-1})\,\circ \,\Delta  \\ 
 \vspace{0.1cm} \\ 
 =\mu \,\circ \,(Id_{i-1} \otimes \mu \otimes Id_{n-i})\circ \lambda _{2n-1}(f_{1}\otimes
f_{2}\otimes \cdots \otimes  f_{2n-1})\circ (Id_{i-1} \otimes \Delta \otimes Id_{n-i})\,\circ \,\Delta ,%
\end{array}%
\end{equation*}%
As $\Delta$ is a $n$-ary totally associative product, we have 
$$
\begin{array}{l}
A(\star )(f_{1}\otimes \cdots \otimes f_{2n-1}) 
\\ = \underset{}{\overset{n-1}{\underset{p=0}{\sum }}}(-1)^{p(n-1)}\mu \circ (Id_{p}
\otimes \mu \otimes Id_{n-1-p}) \circ  \lambda _{2n-1}(f_{1}\otimes \cdots \otimes f_{2n-1})\circ 
(Id_{p}\otimes \Delta \otimes Id_{n-1-p})\circ \Delta  \\ 
\\ = \left( \underset{}{\overset{n-1}{\underset{p=0}{\sum }}}(-1)^{p(n-1)}\mu \circ (Id_{p}
\otimes \mu \otimes Id_{n-1-p}) \right) \circ  \lambda _{2n-1}(f_{1}\otimes \cdots \otimes f_{2n-1})\circ 
(\Delta \otimes Id_{n-1})\circ \Delta=0.  
\end{array}
$$

\section{Some examples of $n$-ary algebras}

\bigskip 1. Let $\mathfrak{g}$ be a Lie algebra. The associator related to
the Lie bracket is 
\begin{equation*}
A(X,Y,Z)=[[X,Y],Z]-[X,[Y,Z]]=[[X,Z],Y].
\end{equation*}%
If $\mathfrak{g}$ is a $4$ step nilpotent Lie algebra the multiplication $%
\mu (X,Y,Z)=A(X,Y,Z)$ is $3$-ary of type $\bullet _{3,3}.$

\bigskip

2. Let $\mu $ a $n$-ary multiplication of type $\bullet _{n,n}$ on a vector
space $V.$ This multiplication is commutative if, for any $v_{i}\in V\,\ ,$%
\begin{equation*}
\sum_{\sigma \in S_{n}}(-1)^{\varepsilon (\sigma )}\mu (v_{\sigma
(1)},\ldots ,v_{\sigma (n)})=0,
\end{equation*}%
where $S_{n}$ is the $n$-order symmetric group and $\varepsilon (\sigma )$
is the signature of the element $\sigma $ of $S_{n}.$ The $3$-ary algebras
of the previous examples are commutative.\ A non-commutative version is
based on the Roby algebras. A Roby algebra is constructed in the following
way: Let $\ V$ be a vector space and $T(V)$ its associated tensor algebra.
For any integer $k$, we consider the ideal $I(V,k)$ of $T(V)$ generated by
the products of symmetric tensors of lenght $k$. The exterior algebra of
order $k$ , or Roby algebra of order $k$, is by definition%
\begin{equation*}
\Lambda (V,k)=T(V)/I(V,k).
\end{equation*}%
For $k=2$ we get the usual exterior algebra. For $k=3,$ the ideal $I(V,3)$
is generated by tensors of type%
\begin{equation*}
\left\{ 
\begin{array}{l}
v_{1}\otimes v_{2}\otimes v_{3}+v_{2}\otimes v_{1}\otimes v_{3}+v_{3}\otimes
v_{2}\otimes v_{1}+v_{1}\otimes v_{3}\otimes v_{2}+v_{2}\otimes v_{3}\otimes
v_{1}+v_{3}\otimes v_{1}\otimes v_{2}, \\ 
v_{1}^{\otimes ^{2}}\otimes v_{2}+v_{2}\otimes v_{1}^{\otimes ^{2}}, \\ 
v_{1}\otimes v_{2}^{\otimes ^{2}}+v_{2}^{\otimes ^{2}}\otimes v_{1},%
\end{array}%
\right.
\end{equation*}%
with distinct vectors $v_{1},v_{2},v_{3}$. If $\mu $ is the multiplication
in $\Lambda (V,3)$, it satisfies%
\begin{equation*}
\left\{ 
\begin{array}{l}
\mu (v_{1},v_{2},v_{3})+\mu (v_{2},v_{1},v_{3})+\mu (v_{3},v_{2},v_{1})+\mu
(v_{1},v_{3},v_{2})+\mu (v_{2},v_{3},v_{1})+\mu (v_{3},v_{1},v_{2})=0 \\ 
\mu (v_{1},v_{1},v_{2})+\mu (v_{2},v_{1},v_{1})=0.%
\end{array}%
\right.
\end{equation*}%
with distinct vector $v_{1},v_{2},v_{3}$. We deduce $\mu
(v_{1},v_{1},v_{1})=0.$ If we now claim that $\mu $ is a multiplication of
type $\bullet _{3,3},$ such algebra is its exterior version.

\bigskip

3. A Poisson algebra of type $\bullet _{3,3}$ can be defined as a
commutative algebra $(V,\mu )$ of type $\bullet _{3,3}$ with a Lie bracket
satisfying%
\begin{equation*}
\lbrack \mu (X,Y,Z),T]=\mu ([X,T],Y,Z)+\mu (X,[Y,T],Z)+\mu (X,Y,[Z,T])
\end{equation*}%
for any $X,Y,Z,T\in V.$ If $V$ is a $\mathbb{Z}_{2}$-graded vector space, we
consider on $V=V_{0}\oplus V_{1}$ a graded Lie bracket which provides $V$
with a super Lie algebra structure.\ Thus this bracket satisfies 
\begin{equation*}
\left\{ 
\begin{array}{lll}
\lbrack X_{1},X_{2}] & = & -[X_{2},X_{1}] \\ 
\left[ X_{1},Y_{2}\right] & = & -[Y_{2},X_{1}] \\ 
\left[ Y_{1},Y_{2}\right] & = & [Y_{2},Y_{1}]%
\end{array}%
\right.
\end{equation*}%
for any $X_{1},X_{2}\in V_{0}$ and $Y_{1},Y_{2}\in V_{1}.$ It also satisfies
the graded Jacobi identity. A superalgebra Poisson structure of type $%
\bullet _{3,3}$ on $V$ $=V_{0}\oplus V_{1}$ is given by a multiplicationt $%
\mu $ of type $\bullet _{3,3}$ and a graded Lie bracket satisfying 
\begin{equation*}
\lbrack \mu (X,Y,Z),T]=\mu ([X,T],Y,Z)+\mu (X,[Y,T],Z)+\mu (X,Y,[Z,T])
\end{equation*}%
An example is given by the $F$-algebras defined in \cite{G.R.T} which are
some generalisation of the superalgebra associated to the super-symmetry. In
fact such algebra (for $F=3)$ is defined on a graded Lie algebra $%
(V=V_{0}\oplus V_{1},[,])$ provided with a commutative multiplication of
type $\bullet _{3,3}$, denoted $\{,,\}$ in this case, and satisfying%
\begin{equation*}
\{V_{i},V_{j},V_{k}\}=0
\end{equation*}%
as soon as $(i,j,k)\neq (1,1,1),$%
\begin{equation*}
\{V_{1},V_{1},V_{1}\}\subseteq V_{0}
\end{equation*}%
and the graded Leibniz relations%
\begin{equation*}
\lbrack
X,\{Y_{1},Y_{2},Y_{3}\}]=\{[X,Y_{1}],Y_{2},Y_{3}\}+\{Y_{1},[X,Y_{2}],Y_{3}%
\}+\{Y_{1},Y_{2},[X,Y_{3}]\}
\end{equation*}%
for any $X\in V_{0}$ et $Y_{1},Y_{2},Y_{3}\in V_{1},$%
\begin{equation*}
\lbrack
Y,\{Y_{1},Y_{2},Y_{3}\}]+[Y_{1},\{Y_{2},Y_{3},Y\}]+[Y_{2},\{Y_{3},Y,Y_{2}%
\}]+[Y_{3},\{Y,Y_{1},Y_{2}\}]=0
\end{equation*}%
for any $Y,$ $Y_{1},Y_{2},Y_{3}\in V_{1}.$

\end{document}